# SENSITIVITY ANALYSIS OF UTILITY-BASED PRICES AND RISK-TOLERANCE WEALTH PROCESSES

BY DMITRY KRAMKOV[1] AND MIHAI SÎRBU

*Carnegie Mellon University and Columbia University*

In the general framework of a semimartingale financial model and a utility function $U$ defined on the positive real line, we compute the first-order expansion of marginal utility-based prices with respect to a "small" number of random endowments. We show that this linear approximation has some important qualitative properties if and only if there is a risk-tolerance wealth process. In particular, they hold true in the following polar cases:

1. for any utility function $U$, if and only if the set of state price densities has a greatest element from the point of view of second-order stochastic dominance;

2. for any financial model, if and only if $U$ is a power utility function ($U$ is an exponential utility function if it is defined on the whole real line).

**1. Introduction.** The valuation of nontraded contingent claims by an economic agent represents a basic problem of financial theory and practice. Given a collection of $m$ such illiquid options, the agent has to determine a "threshold" price $p = (p_i)_{1 \le i \le m}$ such that he/she is willing to buy the $i$th contingent claim at a price less than $p_i$, sell it at a price greater than $p_i$ and do nothing at price $p_i$.

In the framework of a complete financial model, each contingent claim can be replicated by a portfolio of traded securities. Therefore, it admits a uniquely defined arbitrage-free price given as the initial wealth of such a portfolio. While complete financial models have many computational advantages, they are still only idealistic representations (or approximations)

Received April 2005; revised May 2006.

[1]Supported in part by NSF Grants DMS-01-39911 and DMS-05-05414.

of real financial markets as the exact replication of options is usually not possible. Hence, the resulting arbitrage-free prices computed in these models should be used in practice rather cautiously. Indeed, assume for a moment that the illiquid contingent claims can suddenly be bought or sold at a price $p^{\text{trade}}$ which only slightly differs from the price $p$ computed in a complete financial model. The naïve interpretation of the price $p$ leads the investor to take an *infinite* position in the contingent claims (either buy or sell short as many options as possible) which is, clearly, a nonsense from a practical point of view.

The previous discussion motivates the development of pricing rules which, among other factors, depend on the current portfolio $(x, q)$ of the investor $[p = p(x, q)]$, where $x$ is the wealth invested in liquid securities and $q = (q_i)_{1 \leq i \leq m}$ is the vector of his holdings in the nontraded contingent claims. Given such a pricing rule, and assuming that the initial portfolio of the economic agent consists exclusively of liquid wealth $x$ [in other words, it has the form $(x, 0)$], the optimal (static) position $q(x)$ in the illiquid options can be found from the following "equilibrium" equation:

$$(1) \qquad p^{\text{trade}} = p(x - \langle p^{\text{trade}}, q(x) \rangle, q(x)).$$

This equation has a clear economic sense: the agent stops changing his position in the contingent claims as soon as the "model" price reaches the market price. Of course, ideally, the functional dependence $p = p(x, q)$ should be such that (1) has a unique, finite solution.

A formal definition of the prices $p(x, q)$ requires the agent to specify an order relation in the space of portfolio configurations involving random endowments. In other words, given two arbitrary portfolios $(x_i, q_i)$, $i = 1, 2$, he/she should be able to say that $(x_1, q_1)$ is "better" than ("worse" than, "equal" to) $(x_2, q_2)$. The classical approach of financial economics is to define the preferences of the investor with respect to future random payoffs in terms of their expected utilities. In this case, the "quality" of a portfolio $(x, q)$ is expressed as the *maximal expected utility* $u(x, q)$ which can be achieved by investing the liquid amount $x$ in the financial market according to an optimal trading strategy. The pricing rules derived from this approach are commonly called *utility-based* and are currently studied extensively in mathematical finance. The precise definition of the *marginal* utility-based prices $p(x, q)$, the main objects studied in this paper, is given in Definition 1 below.

In most models, explicit computation of marginal utility-based prices can only be performed when the agent holds no contingent claims, that is, for the price $p(x) \triangleq p(x, 0)$; see the papers of Rubinstein [16] and Davis [1]. The case when the agent *does* hold contingent claims is more involved and, except for some rather specific cases (see, e.g., the paper [15] of Musiela and Zariphopoulou), the explicit computations are not possible. This difficulty suggests the performance of an asymptotic analysis of marginal utility-based prices for "small" quantities $q$ of the nontraded contingent claims, that is, the computation of the following *first-order approximation*:

$$(2) \qquad p(x + \Delta x, q) = p(x) + p'(x)\Delta x + D(x)q + o(|\Delta x| + \|q\|).$$

The main goal of this paper is to study this linear approximation in the framework of a general semimartingale financial model and an arbitrary utility function $U$ defined on the positive real line. Note that a similar asymptotic analysis of utility-based prices has been performed by Henderson and Hobson [7], Henderson [6] and Kallsen [12]. The papers [7] and [6] (in a more general case) study the second-order approximation of *total* prices [the quantities $b(x, q)$ and $c(x, q)$ in Appendix A], with respect to the number of claims $q$, in the Black and Scholes model with "basis" risk and for power and exponential utilities. Kallsen [12] studies the first-order approximation of *marginal* prices (similar to our analysis) but in a different and rather nonstandard framework of *local utility maximization*. Furthermore, he proves that in such a model, the first-order correction matrices are symmetric. From a mathematical point of view, our paper is more closely related to [6]. For example, the proof of our key technical Theorem 1 follows a similar "duality based" approach to that used in the proof of Theorem 4.1 in [6].

The principal novelty of our work (to the best of our knowledge) is that in addition to the natural *quantitative* problem of the computation of the vector $p'(x)$ and the matrix $D(x)$ which is accomplished in Theorem 3 below, we are also interested in the following questions of a *qualitative* nature:

1. Are the marginal utility-based prices computed at $q = 0$ locally independent of the initial capital? That is, does

$$(3) \qquad p'(x) = 0$$

hold true?

2. Is the matrix $D(x)$ symmetric, that is, does

$$(4) \qquad D_{ij}(x) = D_{ji}(x), \qquad 1 \le i, j \le m,$$

hold true?

3. Is the matrix $D(x)$ negative semidefinite, that is, does

$$(5) \qquad \langle q, D(x)q \rangle \le 0, \qquad q \in \mathbf{R}^m,$$

hold true?

4. If, in addition, no nonzero linear combination of the contingent claims is replicable (see Assumption 1 below for the exact statement), then does the matrix $D(x)$ have full rank? That is, does

$$(6) \qquad D(x)q = 0 \quad \Longleftrightarrow \quad q = 0, \qquad q \in \mathbf{R}^m,$$

hold true?

In order to motivate these qualitative properties, consider the previously formulated "trading" problem, where the investor can perform a one-time trade in the illiquid options at a price $p^{\text{trade}}$. (For the case of exponential utility, the corresponding problem has been studied recently by Ilhan et al. [10].) As we have already noted, in most cases, it is impossible to compute $p(x, q)$ explicitly and, hence, to provide an exact solution to the "equilibrium" equation (1). If we now approximate the marginal utility-based price $p(x, q)$ by its first-order expansion (2), then (1) becomes

$$(7) \qquad p^{\text{trade}} \approx p(x) - \langle p^{\text{trade}}, q(x) \rangle p'(x) + D(x)q(x).$$

Intuitively, we expect that the optimal (static) position $q(x)$ in the derivative securities should satisfy

$$(8) \qquad \langle q(x), p(x) - p^{\text{trade}} \rangle \geq 0.$$

It can be easily seen that the approximate solution for $q(x)$ implied by (7) is unique and satisfies (8) for any $p^{\text{trade}}$ if and only if $p'(x) = 0$ and $D(x)$ is negative definite, that is, if and only if the qualitative properties (3), (5) and (6) hold true.

In the case of one contingent claim, condition (8) means that $q(x)$ has the same sign as $p(x) - p^{\text{trade}}$. In the case of many contingent claims, the situation is more complicated: some contingent claims $f_i$ might be bought for higher prices than their corresponding threshold prices $p_i(x)$ since the risk coming from buying these claims might be hedged by selling some other derivatives. Therefore, the following question arises: is it possible to find a family of contingent claims $h = (h_i)_{1 \leq i \leq m}$ spanning the same space as the contingent claims $f = (f_i)_{1 \leq i \leq m}$ and such that, for the contingent claims $h$, a change in the traded price of $h_i$ will cause the agent to take a position in the $i$th claim alone? In other words, can we make the corresponding matrix $D(x)$ for $h$ diagonal? This question is directly related to the qualitative property (4).

Perhaps surprisingly, all of the qualitative properties (3)–(6) turn out to be closely related. Moreover, as Theorem 9 below shows, either of these properties holds true if and only if the risk tolerance of optimal terminal wealth (when the investor holds no contingent claims) can be replicated by trading. We call the capital process of such a replicating strategy the *risk-tolerance wealth process*. The latter concept is new and is introduced and studied in Section 5.

The framework of utility-based valuation relies on two basic ingredients: a financial model for traded securities and a utility function $U$ for the preferences of the economic agent. (Strictly speaking, the reference probability measure of the financial model is also a part of the preferences description for the investor; see, e.g., the book of Föllmer and Schied [5], Chapter 2.) In Theorem 6, we show that the risk-tolerance wealth process is well defined [equivalently, the properties (3)–(6) for the marginal utility-based prices hold true] for *any* utility function $U$ if and only if the financial model has the rather remarkable property that the family of state price densities has a greatest element with respect to *second-order stochastic dominance.* Conversely, if we fix the utility function $U$ and require the above properties to hold true for *any* financial model, then, as Theorem 7 shows, $U$ is a power utility (exponential utility, if it is defined on the whole real line). It is interesting to note that most papers on utility-based valuation use a framework that falls into one of these two categories.

The paper is organized as follows. In Section 2, we introduce the financial model and recall the definition of marginal utility-based prices. In Section 3, we perform the *second-order expansion* for the value function $u(x, q)$ of the problem of optimal investment with random endowments. This section is probably the most technically involved and lays the foundation for all subsequent results. The computation of the *sensitivity parameters* $p'(x)$ and $D(x)$ is carried out in Section 4. In Section 5, we define *risk-tolerance wealth processes* and study their existence. Section 6 is dedicated to the study of *qualitative properties* of *marginal utility-based prices*. In order to illustrate the general theory, in Section 7 we work out a specific example of a Black and Scholes model with "basis" risk. For the convenience of future references, in Appendix A we perform the *second-order approximation* for other popular concepts of utility-based prices, namely reservation prices and certainty equivalence values. Finally, Appendix B contains an extension of a well-known result of Jacka [11].

We conclude this section by introducing some notation to be used throughout the paper. For vectors $a = (a_i)_{1 \leq i \leq m}$ and $b = (b_i)_{1 \leq i \leq m}$ and matrix $C = (C_{ij})_{1 \leq i, j \leq m}$, we denote

$$\langle a, b \rangle \triangleq \sum_{i=1}^{m} a_i b_i, \qquad \|a\| \triangleq \sqrt{\langle a, a \rangle}, \qquad \|C\| \triangleq \sup_{a \in \mathbf{R}^m} \frac{\|Ca\|}{\|a\|}.$$

Given a filtered probability space $(\Omega, \mathcal{F}, (\mathcal{F}_t)_{0 \leq t \leq T}, \mathbb{P})$ and a number $p \geq 1$, we denote by $\mathbf{H}^p(\mathbb{P})$ the space of uniformly integrable martingales under $\mathbb{P}$ with the norm $(\mathbb{E}|M_T^*|^p)^{1/p}$, where $M_T^* \triangleq \sup_{0 \leq t \leq T} |M_t|$. We also denote

$$\mathbf{H}_0^p(\mathbb{P}) \triangleq \{M \in \mathbf{H}^p(\mathbb{P}) : M_0 = 0\}.$$

**2. The definitions of basic concepts.** We study a financial framework similar to the one in [8], referring to this paper for more details and references. We consider a model of a security market which consists of $d + 1$ assets, one bond and $d$ stocks. We work in discounted terms, that is, we suppose that the price of the bond is constant and denote by $S = (S^i)_{1 \leq i \leq d}$ the price process of the $d$ stocks. The process $S$ is assumed to be a semi-martingale on a filtered probability space $(\Omega, \mathcal{F}, (\mathcal{F}_t)_{0 \le t \le T}, \mathbb{P})$. Here $T$ is a finite time horizon and $\mathcal{F} = \mathcal{F}_T$.

A (self-financing) portfolio is defined as a pair $(x, H)$, where the constant $x$ represents the initial capital and $H = (H^i)_{1 \le i \le d}$ is a predictable $S$-integrable process, $H_t^i$ specifying how many units of asset $i$ are held in the portfolio at time $t$. The wealth process $X = (X_t)_{0 \le t \le T}$ of the portfolio evolves in time as the stochastic integral of $H$ with respect to $S$:

$$(9) \qquad X_t = x + \int_0^t H_u \, dS_u, \qquad 0 \le t \le T.$$

We denote by $\mathcal{X}(x)$ the family of wealth processes with nonnegative capital at any instant and with initial value equal to $x$:

$$(10) \qquad \mathcal{X}(x) \triangleq \{X \ge 0 \colon X \text{ is defined by (9)}\}.$$

A nonnegative wealth process is said to be *maximal* if its terminal value cannot be dominated by that of any other nonnegative wealth process with the same initial value. In general, a wealth process $X$ is called *maximal* if it admits a representation of the form

$$X = X' - X'',$$

where both $X'$ and $X''$ are nonnegative maximal wealth processes. A wealth process $X$ is said to be *acceptable* if it admits a representation as above, where both $X'$ and $X''$ are nonnegative wealth processes and, in addition, $X''$ is maximal. For a detailed discussion of maximal and acceptable processes, we refer the reader to [3].

A probability measure $\mathbb{Q} \sim \mathbb{P}$ is called an *equivalent local martingale measure* if any $X \in \mathcal{X}(1)$ is a local martingale under $\mathbb{Q}$. The family of equivalent local martingale measures is denoted by $\mathcal{Q}$. We assume throughout that

$$(11) \qquad \qquad \qquad \mathcal{Q} \ne \varnothing.$$

This rather mild condition is equivalent to the absence of arbitrage opportunities in the model; see [2] and [4] for precise statements and further references. In particular, (11) implies that a constant positive process is maximal.

In addition to the set of traded securities, we consider a family of $m$ non-traded European contingent claims with payment functions $f = (f_i)_{1 \le i \le m}$, which are $\mathcal{F}$-measurable random variables, and maturity $T$. We assume that this family is dominated by the terminal value of some nonnegative wealth process $X$, that is, that

$$(12) \qquad \qquad \|f\| \triangleq \sqrt{\sum_{i=1}^{m} f_i^2} \le X_T.$$

It is well known (see [2], Theorem 5.7) that (12) is equivalent to the following integrability condition:

$$(13) \qquad \sup_{\mathbb{Q} \in \mathcal{Q}} \mathbb{E}_{\mathbb{Q}}[\|f\|] < \infty.$$

Recall that a random variable $g$ is *replicable* if there exists a maximal wealth process $X$ such that $X_T = g$. Provided that it exists, such a process $X$ is unique and is called the *replication process* for $g$. Hereafter, we shall often impose (without loss of generality) the following assumption on the family of random endowments $f$:

ASSUMPTION 1. For any $q \in \mathbf{R}^m$ such that $q \neq 0$, the random variable $\langle q, f \rangle \triangleq \sum_{i=1}^m q_i f_i$ is not replicable.

For $x \in \mathbf{R}$ and $q \in \mathbf{R}^m$, we denote by $\mathcal{X}(x, q)$ the set of acceptable processes with initial capital $x$ whose terminal values dominate $-\langle q, f \rangle$, that is,

$$\mathcal{X}(x, q) \triangleq \{X : X \text{ is acceptable, } X_0 = x \text{ and } X_T + \langle q, f \rangle \geq 0\}.$$

From the definition of acceptable processes, we deduce that

$$\mathcal{X}(x, 0) = \mathcal{X}(x),$$

where the set $\mathcal{X}(x)$ is defined in (10). The set of points $(x, q)$ where $\mathcal{X}(x, q)$ is not empty is a closed convex cone in $\mathbf{R}^{m+1}$. We denote by $\mathcal{K}$ the interior of this cone, that is,

$$\mathcal{K} \triangleq \operatorname{int}\{(x, q) \in \mathbf{R}^{m+1} : \mathcal{X}(x, q) \neq \varnothing\}.$$

Note (see [8], Lemma 1), that (12) is equivalent to the fact that $\mathcal{K}$ contains any point $(x, q)$ such that $x > 0$ and $q = 0$:

$$(14) \qquad (x, 0) \in \mathcal{K}, \qquad x > 0.$$

We also consider an economic agent in our model, whose preferences over terminal wealth are modeled by a utility function $U : (0, \infty) \to (-\infty, \infty)$. The function $U$ is assumed to be strictly concave, strictly increasing, continuously differentiable and to satisfy the Inada conditions:

$$(15) \qquad U'(0) = \lim_{x \to 0} U'(x) = \infty, \qquad U'(\infty) = \lim_{x \to \infty} U'(x) = 0.$$

Following [13], we assume that the asymptotic elasticity of $U$ is strictly less than 1, that is, that

$$(16) \qquad \limsup_{x \to \infty} \frac{x U'(x)}{U(x)} < 1.$$

Hereafter, we set $U(0) = \lim_{x \to 0} U(x)$ and $U(x) = -\infty$ for all $x < 0$.

Assume that the agent has some initial capital $x$ and quantities $q = (q_i)_{1 \leq i \leq m}$ of the contingent claims $f$ such that $(x, q) \in \mathcal{K}$. The quantities $q$ of the contingent claims will be held constant up to the maturity. On the contrary, the capital $x$ can be freely invested in the stocks and the bond according to some dynamic strategy. The maximal expected utility that the agent can achieve by trading in the financial market is thus given by

$$(17) \qquad u(x, q) \triangleq \sup_{X \in \mathcal{X}(x, q)} \mathbb{E}[U(X_T + \langle q, f \rangle)], \qquad (x, q) \in \mathcal{K}.$$

We shall use a shorter notation $u(x)$ for the value function in the case without random endowments, that is,

$$(18) \qquad u(x) \triangleq u(x, 0) = \sup_{X \in \mathcal{X}(x)} \mathbb{E}[U(X_T)], \qquad x > 0.$$

To exclude the trivial case, we shall assume that

$$(19) \qquad u(x) < \infty \qquad \text{for some } x > 0.$$

Since $u = u(x, q)$ is a concave function on the open set $\mathcal{K}$ and as, by (14), any point $(x, 0)$ belongs to $\mathcal{K}$, this condition implies that

$$(20) \qquad u(x, q) < \infty \qquad \text{for all } (x, q) \in \mathcal{K}.$$

An important role in the analysis to follow will be played by the *marginal utility of the terminal wealth of the optimal investment strategy*, that is, by the random variable $U'(X_T(x))$, where $X(x)$ is the solution to (18). In practice, it is often easier to compute $U'(X_T(x))$ than to compute the terminal wealth $X_T(x)$ itself. This computation is based on the fact that the former represents the solution of the dual problem to (18). More precisely,

$$(21) \qquad U'(X_T(x)) = Y_T(u'(x)),$$

where $Y(y)$ is the solution to the following optimization problem:

$$(22) \qquad v(y) \triangleq \inf_{Y \in \mathcal{Y}(y)} \mathbb{E}[V(Y_T)], \qquad y > 0.$$

Here $V$ is the convex conjugate to $U$, that is,

$$V(y) \triangleq \sup_{x > 0} \{ U(x) - xy \}, \qquad y > 0$$

and $\mathcal{Y}(y)$ is the family of nonnegative supermartingales $Y$ such that $Y_0 = y$ and $XY$ is a supermartingale for all $X \in \mathcal{X}(1)$. Note that $\mathcal{Y}(1)$ contains the density processes of all $\mathbb{Q} \in \mathcal{Q}$.

According to [13], Theorem 2.2, conditions (11), (16) and (19) imply that the value functions $u(x)$ and $v(y)$ are conjugate, that is,

$$(23) \qquad v(y) = \sup_{x > 0} \{ u(x) - xy \}, \qquad y > 0,$$

and that they are continuously differentiable on $(0, \infty)$. In addition, the solutions $X(x)$ and $Y(y)$ to (18) and (22) exist and are unique for all $x > 0$ and $y > 0$ and, if $y = u'(x)$, then the product $X(x)Y(y)$ is a martingale. Hereafter, we shall use these results without further comment.

We now introduce the primary object of our study:

DEFINITION 1. Let $(x, q) \in \mathcal{K}$ be the initial portfolio of the agent. A vector $p \in \mathbf{R}^m$ is called a *marginal utility-based price* (for the contingent claims $f$) at $(x, q)$ if

$$u(x, q) \geq u(x', q')$$

for any $(x', q') \in \mathcal{K}$ such that

$$x + \langle q, p \rangle = x' + \langle q', p \rangle.$$

The interpretation of this definition is that the agent's holdings $q$ in $f$ are optimal in the model where the contingent claims can be traded at the marginal utility-based price $p$ at time zero. Marginal utility-based prices are, in general, not unique. We denote by $\mathcal{P}(x, q)$ the set of marginal utility-based prices at $(x, q) \in \mathcal{K}$. Using standard arguments from the theory of convex functions, we deduce that

$$(24) \qquad p \in \mathcal{P}(x, q) \quad \Longleftrightarrow \quad p = \frac{r}{y} \qquad \text{for some } (y, r) \in \partial u(x, q),$$

where $\partial u(x, q)$ is the subdifferential of the value function $u(x, q)$ defined in (17) at $(x, q)$. In particular, the set $\mathcal{P}(x, q)$ is a singleton [i.e., marginal utility-based prices are uniquely defined at $(x, q)$] if and only if $u(x, q)$ is differentiable at this point.

Explicit formulas for marginal utility-based prices are usually available only when $q = 0$. It was shown in [9] that the set $\mathcal{P}(x, 0)$ is a singleton and that its unique element, denoted by $p(x)$, has the representation

$$(25) \qquad p(x) = \frac{1}{y} \mathbb{E}[Y_T(y)f], \qquad y = u'(x),$$

provided that there is a positive wealth process $X$ that dominates $f$ according to (12) and such that

$$(26) \qquad XY(y) \text{ is a martingale on } [0, T], \qquad y = u'(x).$$

The main goal of this paper is to study the *first-order approximation* of the marginal utility-based prices in the neighborhood of $(x, 0)$,

$$(27) \qquad p = p(x) + p'(x)\Delta x + D(x)q + o(|\Delta x| + \|q\|), \qquad p \in \mathcal{P}(x + \Delta x, q),$$

or, more precisely, to find $p'(x)$ and $D(x)$ such that

$$(28) \qquad \lim_{|\Delta x| + \|q\| \to 0} \sup_{p \in \mathcal{P}(x + \Delta x, q)} \frac{\|p - (p(x) + p'(x)\Delta x + D(x)q)\|}{|\Delta x| + \|q\|} = 0.$$

Here, $p'(x) = (p_i'(x))_{1 \leq i \leq m}$ is an $m$-dimensional vector and $D(x) = (D_{ij}(x))_{1 \leq i,j \leq m}$ is an $m \times m$ matrix. Note that if the marginal utility-based prices are uniquely defined in the vicinity of $(x, 0)$, then $p'(x)$ and $D(x)$ are their first derivatives (in the classical sense) with respect to $x$ and $q$ computed at this point.

## 3. Second-order approximation for $u(x, q)$ in the vicinity of $(x, 0)$.
Throughout this section we fix the initial capital $x > 0$. Our study of the approximation (27) will be based on the second-order expansion of the value function $u(x, q)$ defined in (17) in the neighborhood of $(x, 0)$,

$$(29) \qquad \begin{aligned} u(x + \Delta x, q) = u(x) &+ u'(x)\Delta x + \langle u_q(x, 0), q \rangle \\ &+ \tfrac{1}{2} \begin{pmatrix} \Delta x & q \end{pmatrix} G(x) \begin{pmatrix} \Delta x \\ q \end{pmatrix} + o((\Delta x)^2 + \|q\|^2). \end{aligned}$$

Here, as the notation suggests, $u'(x)$ is the first derivative of the value function $u(x)$ defined in (18) and $u_q(x, 0)$ is the vector of the partial derivatives of $u(x, q)$ with respect to $q$ computed at $(x, 0)$. Recall that representations (24) and (25) for the marginal utility-based prices $p(x)$ at $(x, 0)$ imply that

$$(30) \qquad u_q(x, 0) = yp(x) = \mathbb{E}[Y_T(y)f], \qquad y = u'(x),$$

where $Y(y)$ is the solution to (22). Further, in this approximation, $G(x) = (G_{ij}(x))_{0 \leq i,j \leq m}$ is the symmetric negative semidefinite matrix of the second-order derivatives of $u(x, q)$ at $(x, 0)$ in the sense that

$$(31) \qquad \lim_{|\Delta x| + \|q\| \to 0} \sup_{z \in \partial u(x + \Delta x, q)} \frac{\|z - \binom{u'(x)}{u_q(x,0)} - G(x)\binom{\Delta x}{q}\|}{|\Delta x| + \|q\|} = 0.$$

[A short proof of the implication (29) $\Longrightarrow$ (31) is given in Lemma 8 below.]

To state the results, we need to impose additional conditions on the utility function $U$, the price process $S$ and the random endowments $f$. The first two assumptions were previously formulated in [14] and are essential for the two-times differentiability of the value functions $u(x)$ and $v(y)$ defined in (18) and (22). We refer to this paper for more details.

ASSUMPTION 2. The utility function $U$ is two times continuously differentiable on $(0, \infty)$ and its relative risk-aversion coefficient

$$(32) \qquad A(x) \triangleq -\frac{xU''(x)}{U'(x)}, \qquad x > 0,$$

is uniformly bounded away from zero and infinity, that is, there exist constants $c_1 > 0$ and $c_2 < \infty$ such that

$$(33) \qquad c_1 < A(x) < c_2, \qquad x > 0.$$

It is easy to show that this assumption implies both the Inada conditions (15) and condition (16) on the asymptotic elasticity; see [14], Lemma 3. Note that Assumption 2 is equivalent to the following symmetric condition on the conjugate function $V$:

$$(34) \qquad \frac{1}{c_2} < B(y) < \frac{1}{c_1}, \qquad y > 0,$$

where

$$(35) \qquad B(y) \triangleq -\frac{yV''(y)}{V'(y)}, \qquad y > 0.$$

Note also that

$$(36) \qquad B(U'(x)) = \frac{1}{A(x)}$$

is the relative *risk-tolerance* coefficient of $U$ computed at $x$.

Following [14], we call a $d$-dimensional semimartingale $R$ *sigma-bounded* if there is a strictly positive, predictable (one-dimensional) process $h$ such that the stochastic integral $\int h\,dR$ is well defined and locally bounded.

ASSUMPTION 3. The price process of the traded securities discounted by the solution $X(x)$ to (18), that is the $(d+1)$-dimensional semimartingale

$$(37) \qquad S^{X(x)} \triangleq \left(\frac{1}{X(x)}, \frac{S}{X(x)}\right),$$

is sigma-bounded.

We refer to [14], Theorem 3 for sufficient conditions that ensure the validity of this assumption. In particular, this assumption is satisfied if $S$ is a continuous process or if the original (incomplete) model can be extended to a complete one by adding a finite number of securities.

To facilitate the formulation of the assumption on the random endowments $f$, we change the numéraire from the bond to the normalized optimal wealth process $X(x)/x$ and denote by

$$(38) \qquad g_i(x) \triangleq x\frac{f_i}{X_T(x)}, \qquad 1 \le i \le m,$$

the payoffs of the European options discounted by $X(x)/x$. To simplify future notation, we add to this family the additional component

$$(39) \qquad g_0(x) \triangleq 1.$$

Let $\mathbb{R}(x)$ be the probability measure on $(\Omega, \mathcal{F})$ whose Radon–Nikodym derivative under $\mathbb{P}$ is given by

$$(40) \qquad \frac{d\mathbb{R}(x)}{d\mathbb{P}} = \frac{X_T(x)Y_T(y)}{xy}, \qquad y = u'(x).$$

Recalling that $\mathbf{H}_0^2(\mathbb{R}(x))$ stands for the space of square integrable martingales under $\mathbb{R}(x)$ with initial value 0, we denote

$$(41) \qquad \mathcal{M}^2(x) \triangleq \left\{ M \in \mathbf{H}_0^2(\mathbb{R}(x)) \colon M = \int H\,dS^{X(x)} \right\},$$

where $S^{X(x)}$ was defined in (37). Note that if $M \in \mathcal{M}^2(x)$, then $X(x)M$ is a wealth process under the original numéraire, that is, a stochastic integral with respect to $S$.

ASSUMPTION 4. There exists a constant $c > 0$ and a process $M \in \mathcal{M}^2(x)$ such that

$$(42) \qquad \sum_{i=1}^m |g_i(x)| \le c + M_T.$$

REMARK 1. Assumption 4 implies that both conditions (12) and (26) on the random endowments $f$ are satisfied. Indeed, these conditions together take the form (42), where $M$ is a stochastic integral with respect to $S^{X(x)}$ and a uniformly integrable martingale under $\mathbb{R}(x)$.

The matrix $G(x)$ in (29) will be computed in terms of the appropriately defined derivatives of the solutions $X(x,q)$ to (17). To characterize these derivative processes, we consider the following optimization problems:

$$(43) \qquad a_i(x) \triangleq \inf_{M \in \mathcal{M}^2(x)} \mathbb{E}_{\mathbb{R}(x)}[A(X_T(x))(g_i(x) + M_T)^2], \qquad 0 \le i \le m,$$

where the function $A$ is defined in (32). It is easy to show that under Assumptions 2 and 4, the solutions $M^i(x)$ to (43), $0 \le i \le m$, exist and are unique. Using this notation, we set

$$(44) \qquad Z^i(x) \triangleq \frac{X(x)}{x} M^i(x), \qquad 0 \le i \le m,$$

$$(45) \qquad X'(x) \triangleq \frac{X(x)}{x}(1 + M^0(x)) = \frac{X(x)}{x} + Z^0(x).$$

From the definition of the set $\mathcal{M}^2(x)$, we deduce that $Z^i(x)$ and $X'(x)$ are wealth processes (i.e., stochastic integrals with respect to $S$) and that if $y = u'(x)$, then $Z^i(x)Y(y)$ and $X'(x)Y(y)$ are uniformly integrable martingales under $\mathbb{P}$.

The next theorem is the main result of this section.

THEOREM 1. *Let $x > 0$. Assume* (11) *and* (19) *and also that Assumptions* 2, 3 *and* 4 *hold true. Then:*

1. *The terminal values of the wealth processes $X'(x)$ and $(Z^i(x))_{1 \leq i \leq m}$ defined in* (45) *and* (44) *are the partial derivatives of the terminal values of the solutions $X(x, q)$ to* (17) *computed at $(x, 0)$, that is*

$$
\begin{aligned}
(46) \quad \lim_{|\Delta x| + \|q\| \to 0} \frac{1}{|\Delta x| + \|q\|} \bigg| &X_T(x + \Delta x, q) \\
&- X_T(x) - X'_T(x)\Delta x - \sum_{i=1}^{m} Z^i_T(x) q_i \bigg| = 0,
\end{aligned}
$$

*where the convergence takes place in $\mathbb{P}$-probability.*

2. *The value function $u(x, q)$ defined in* (17) *has the second-order expansion* (29)*, where $u_q(x, 0)$ is given by* (30) *and, for $0 \leq i, j \leq m$,*

$$
\begin{aligned}
(47) \quad G_{ij}(x) &= \mathbb{E}[U''(X_T(x))(Z^i_T(x) + f_i)(Z^j_T(x) + f_j)] \\
&= -\frac{u'(x)}{x} \mathbb{E}_{\mathbb{R}(x)}[A(X_T(x))(M^i_T(x) + g_i(x))(M^j_T(x) + g_j(x))],
\end{aligned}
$$

*where $(M^i(x))_{0 \leq i \leq m}$ are the solutions to* (43) *and where we have denoted $f_0 \triangleq X_T(x)/x$ so that $Z^0_T(x) + f_0 = X'_T(x)$.*

The second-order expansion of $u(x, q)$ will be studied *simultaneously* with the second-order expansion of its *dual* function $v(y, r)$ given by

$$
(48) \qquad v(y, r) \triangleq \inf_{Y \in \mathcal{Y}(y, r)} \mathbb{E}[V(Y_T)], \qquad (y, r) \in \mathcal{L}.
$$

Here, the *dual domain* $\mathcal{L}$ is defined as the relative interior of the polar cone of $-\mathcal{K}$,

$$
\mathcal{L} \triangleq \mathrm{ri}\{(y, r) \in \mathbf{R}^{m+1} : xy + \langle q, r \rangle \geq 0 \text{ for all } (x, q) \in \mathcal{K}\}
$$

and $\mathcal{Y}(y, r)$ [for $(y, r) \in \mathcal{L}$] is the set of nonnegative supermartingales $Y \in \mathcal{Y}(y)$ such that for any $(x, q) \in \mathcal{K}$ and any $X \in \mathcal{X}(x, q)$, we have

$$
\mathbb{E}[Y_T(X_T + \langle q, f \rangle)] \leq xy + \langle q, r \rangle.
$$

According to [8], Theorem 1, the functions $u(x, q)$ and $v(y, r)$ are conjugate:

$$
\begin{aligned}
(49) \quad u(x, q) &= \inf_{(y, r) \in \mathcal{L}} \{v(y, r) + xy + \langle q, r \rangle\}, \qquad (x, q) \in \mathcal{K}, \\
v(y, r) &= \sup_{(x, q) \in \mathcal{K}} \{u(x, q) - xy - \langle q, r \rangle\}, \qquad (y, r) \in \mathcal{L}.
\end{aligned}
$$

To avoid trivial complications, we assume in the study of (48) that Assumption 1 holds true. In this case (see [8], Lemmas 3 and 7) $\mathcal{L}$ is an open cone in $\mathbf{R}^{m+1}$ and $v(y,r)$ is continuously differentiable on $\mathcal{L}$.

The second-order expansion for $v(y,r)$ will be performed at $(y, r(y))$, where, for $y > 0$,

$$(50) \qquad r(y) = \mathbb{E}[Y_T(y)f].$$

Taking into account (30), we deduce that for $y = u'(x)$,

$$(51) \qquad (y, r(y)) = \partial u(x, 0)$$

and hence, by (49) and (23),

$$v(y, r(y)) = u(x, 0) - xy - \langle 0, r(y) \rangle = u(x) - xy = v(y).$$

It follows that

$$Y(y) = Y(y, r(y)),$$

where $Y(y)$ is the solution to (22) and $Y(y, r)$ is the solution to (48).

In Theorem 2, we shall show that

$$v(y + \Delta y, r(y) + \Delta r) = v(y) + v'(y)\Delta y + \tfrac{1}{2} \begin{pmatrix} \Delta y & \Delta r \end{pmatrix} H(y) \begin{pmatrix} \Delta y \\ \Delta r \end{pmatrix}$$
$$(52) \qquad\qquad\qquad + o((\Delta y)^2 + \|\Delta r\|^2),$$

where $v'(y)$ is the first derivative of $v(y)$ and $H(y)$ is the matrix of the second derivatives of $v(y, r)$ at $(y, r(y))$. In Theorem 2, the matrix $H(y)$ will be computed in terms of the derivatives of $Y(y, r)$ at $(y, r(y))$.

For $y = u'(x)$, we denote by $\mathcal{N}^2(y)$ the orthogonal complement to $\mathcal{M}^2(x)$ in $\mathbf{H}_0^2(\mathbb{R}(x))$, that is,

$$(53) \qquad \mathcal{N}^2(y) \triangleq \{N \in \mathbf{H}_0^2(\mathbb{R}(x)) : NM \text{ is a } \mathbb{R}(x)\text{-martingale}$$
$$\text{for all } M \in \mathcal{M}^2(x)\},$$

and define the following *affine* subspaces of the space $\mathbf{H}^2(\mathbb{R}(x))$ of square integrable martingales under $\mathbb{R}(x)$:

$$\mathcal{N}_0^2(y) \triangleq \{N \in 1 + \mathcal{N}^2(y) : \mathbb{E}_{\mathbb{R}(x)}[N_T g_j(x)] = 0, 1 \le j \le m\},$$

$$\mathcal{N}_i^2(y) \triangleq \{N \in \mathcal{N}^2(y) : \mathbb{E}_{\mathbb{R}(x)}[N_T g_j(x)] = \delta_{ij}, 1 \le j \le m\}, \qquad 1 \le i \le m,$$

where $\delta_{ij}$ is the standard notation for $I_{\{i=j\}}$. In Lemma 4, we shall show that under the conditions of Theorem 1 and Assumption 1, the sets $\mathcal{N}_i^2(y)$ are nonempty. Hence, we can introduce the following optimization problem:

$$(54) \qquad b_i(y) \triangleq \inf_{N \in \mathcal{N}_i^2(y)} \mathbb{E}_{\mathbb{R}(x)}[B(Y_T(y))(N_T)^2], \qquad 0 \le i \le m,$$

where $B$ is given by (35). As we shall see in Lemma 7, this problem is, in fact, dual to the one defined in (43). Taking into account (34), it is easy to deduce that the solutions $N^i(y)$ to (54) exist and are unique. Using this notation, we define

$$(55) \qquad W^i(y) \triangleq \frac{Y(y)}{y} N^i(y), \qquad 0 \leq i \leq m.$$

We remark that the products $X(x)W^i(y)$ and $X'(x)W^i(y)$ are uniformly integrable martingales under $\mathbb{P}$.

THEOREM 2. *Assume that the conditions of Theorem 1 and Assumption 1 hold true. Let $x > 0$ and $(y, r(y))$ be given by (51). Then:*

1. *The terminal values of the $(m+1)$-dimensional process $W(y)$ defined in (55) are the partial derivatives at $(y, r(y))$ of the terminal values of the solutions $Y(y, r)$ to (48),*

$$
\begin{aligned}
(56) \quad \lim_{|\Delta y| + \|\Delta r\| \to 0} \frac{1}{|\Delta y| + \|\Delta r\|} \bigg| & Y_T(y + \Delta y, r(y) + \Delta r) \\
& - Y_T(y) - W_T^0(y)\Delta y - \sum_{i=1}^m W_T^i(y)\Delta r_i \bigg| = 0,
\end{aligned}
$$

*where the above limit takes place in $\mathbb{P}$-probability.*

2. *The value function $v(y, r)$ defined in (48) has the second-order expansion (52), where for $0 \leq i, j \leq m$,*

$$
\begin{aligned}
(57) \quad H_{ij}(y) &= \mathbb{E}[V''(Y(y))W_T^i(y)W_T^j(y)] \\
&= \frac{x}{y}\mathbb{E}_{\mathbb{R}(x)}[B(Y_T(y))N_T^i(y)N_T^j(y)]
\end{aligned}
$$

*and $N^i(y)$ is the solution to (54). The matrix $H(y)$ is positive definite and its inverse is given by*

$$(58) \qquad (H(y))^{-1} = -G(x),$$

*where $G(x)$ was defined in (47).*

3. *The derivative processes $X'(x)$ and $Z(x) \triangleq (Z^i(x))_{1 \leq i \leq m}$ defined in (45) and (44) and $W(y) = (W^i(y))_{0 \leq i \leq m}$ are related by (in vector notation)*

$$(59) \qquad U''(X_T(x)) \begin{pmatrix} X_T'(x) \\ Z_T(x) + f \end{pmatrix} = G(x)W_T(y).$$

The proofs of Theorems 1 and 2 will be given simultaneously and will follow from a series of lemmas. The basic idea is similar to the approach used in [6]. Namely, we approximate the primal value function to the second-order "from below" [equation (67)] and approximate up to second-order the dual value function "from above" [equation (68)]. If the two approximations "match" (i.e., the Hessian matrices satisfy the usual relation corresponding to convex conjugate functions), then we have a "two-sided" second-order approximation for both value functions. The second-order approximation of value functions is directly related to the *first-order approximation* of optimal investment strategies and dual optimizers. We compute these first-order approximations as the result of some auxiliary optimization problems, namely (43), (44), (54) and (55).

As before, we fix $x > 0$ and denote $y = u'(x)$. We start with a technical result which is a slight improvement over Lemma 4 in [14]:

LEMMA 1. *Assume that the utility function $U = U(x)$ satisfies Assumption* 2. *Let $\xi > 0$ be a random variable such that*

$$\mathbb{E}[|U(\xi)|] < \infty \tag{60}$$

*and let $\eta = (\eta^i)_{0 \leq i \leq m}$ be a $(m+1)$-dimensional random vector such that*

$$\eta^i \geq -K\xi, \qquad 0 \leq i \leq m,$$

*for some $K > 0$ and*

$$\mathbb{E}[-U''(\xi)\|\eta\|^2] < \infty. \tag{61}$$

*Then the function*

$$w(s) \triangleq \mathbb{E}[U(\xi + \langle s, \eta \rangle)]$$

*is well defined for $s \in \mathbf{R}^{m+1}$ such that*

$$0 \leq s^i < \frac{1}{K(m+1)}, \qquad 0 \leq i \leq m, \tag{62}$$

*and*

$$w(s) = w(0) + \mathbb{E}[U'(\xi)\langle \eta, s \rangle] + \tfrac{1}{2}\mathbb{E}[U''(\xi)|\langle \eta, s \rangle|^2] + o(\|s\|^2).$$

PROOF. As in the proof of [14], Lemma 4, we deduce from (60) and Assumption 2 that

$$\mathbb{E}[\xi U'(\xi)] < \infty.$$

Defining the measure $\mathbb{R}(\xi)$ by

$$\frac{d\mathbb{R}(\xi)}{d\mathbb{P}} = \frac{\xi U'(\xi)}{\mathbb{E}[\xi U'(\xi)]},$$

we observe, taking (33) into account, that condition (61) means, in fact, that $\eta/\xi \in \mathbf{L}^2(\mathbb{R}(\xi))$. This implies $\eta/\xi \in \mathbf{L}^1(\mathbb{R}(\xi))$, which means

$$(63) \qquad \mathbb{E}[U'(\xi)\|\eta\|] < \infty.$$

Let $0 < b < \frac{1}{K(m+1)}$ and $\delta \triangleq 1 - K(m+1)b > 0$. If

$$(64) \qquad 0 \leq s^i \leq b, \qquad 0 \leq i \leq m,$$

then

$$\xi + \langle s, \eta \rangle \geq \xi - K\left(\sum_{i=0}^m s^i\right)\xi \geq \delta\xi.$$

For $s, t$ satisfying (64), we obtain

$$
\begin{aligned}
(65) \qquad \frac{|U(\xi + \langle s, \eta \rangle) - U(\xi + \langle t, \eta \rangle)|}{\|s - t\|} &= \frac{U'(\zeta_{s,t})|\langle s - t, \eta \rangle|}{\|s - t\|} \\
&\leq U'(\delta\xi)\|\eta\| \\
&\leq C U'(\xi)\|\eta\|,
\end{aligned}
$$

where $\zeta_{s,t}$ takes intermediate values between $\xi + \langle s, \eta \rangle$ and $\xi + \langle t, \eta \rangle$ (therefore, $\zeta_{s,t} \geq \delta\xi$) and, at the last step, we used the inequality $U'(\delta\xi) \leq CU'(\xi)$ that follows from (33) by [14], Lemma 3. Similarly, from the last inequality and Assumption 2, we deduce for $s, t$ satisfying (64) that

$$
\begin{aligned}
(66) \qquad \frac{\|U'(\xi + \langle s, \eta \rangle)\eta - U'(\xi + \langle t, \eta \rangle)\eta\|}{\|s - t\|} &\leq |U''(\zeta_{s,t})|\|\eta\|^2 \\
&\leq C_1 \frac{U'(\zeta_{s,t})}{\zeta_{s,t}}\|\eta\|^2 \\
&\leq C_2 \frac{U'(\xi)}{\xi}\|\eta\|^2 \\
&\leq C_3|U''(\xi)|\|\eta\|^2,
\end{aligned}
$$

where $(C_i)_{1 \leq i \leq 3}$ are some positive constants. Taking into account (65) and (66) together with (63) and (61), we can now apply the dominated convergence theorem to complete the proof. $\quad\square$

Using the fact that $u(x, q)$ is the value function of a *maximization* problem, we shall deduce next the *lower* bound,

$$
\begin{aligned}
(67) \qquad u(x + \Delta x, q) &\geq u(x) + y\Delta x + \langle r(y), q \rangle \\
&\quad + \tfrac{1}{2}\begin{pmatrix} \Delta x & q \end{pmatrix} G(x) \begin{pmatrix} \Delta x \\ q \end{pmatrix} + o((\Delta x)^2 + \|q\|^2),
\end{aligned}
$$

where $y = u'(x)$, the vector $r(y)$ is given by (50) and the matrix $G(x)$ is defined by (47). We start with a preliminary result:

LEMMA 2. *Assume that the conditions of Theorem* 1 *hold true. Let* $M \in \mathcal{M}^2(x)$ *and let* $h$ *be a random variable such that*

$$|h| \leq c + N_T,$$

*for some constant* $c > 0$ *and a process* $N \in \mathcal{M}^2(x)$.

*Then there is a sequence* $(M^n)_{n \geq 1}$ *in* $\mathcal{M}^2(x)$ *that converges to* $M$ *in* $\mathbf{H}_0^2(\mathbb{R}(x))$ *and such that for any* $n \geq 1$, *the random variable* $M_T^n + h$ *is bounded below.*

PROOF. Since the stock process $S^{X(x)}$ [with respect to the new numéraire $X(x)/x$] is *sigma-bounded*, there is a sequence $(L^n)_{n \geq 1}$ in $\mathcal{M}^2(x)$ that converges to $M - N$ in $\mathbf{H}_0^2(\mathbb{R}(x))$ and such that any $L^n$ is bounded. We can now choose $M^n \triangleq N + L^n$ to complete the proof. □

LEMMA 3. *Assume that the conditions of Theorem* 1 *are satisfied. Then* (67) *holds true.*

PROOF. We give the proof of (67) for $(\Delta x, q)$ in the first orthant of $\mathbf{R}^{m+1}$. The other cases are treated similarly.

For $0 \leq i \leq m$, let $M^i$ be an element of $\mathcal{M}^2(x)$ such that $M_T^i + g_i(x)$ is bounded below. It is easy to see that for sufficiently small $(\Delta x, q)$ in the first orthant of $\mathbf{R}^{m+1}$,

$$X \triangleq \frac{X(x)}{x}\left(x + \Delta x(1 + M^0) + \sum_{i=1}^{m} q_i M^i\right) \in \mathcal{X}(x + \Delta x, q)$$

and hence,

$$u(x + \Delta x, q) \geq \mathbb{E}[U(X_T + \langle q, f \rangle)]$$
$$= \mathbb{E}[U(X_T(x) + (X_T - X_T(x)) + \langle q, f \rangle)]$$
$$= \mathbb{E}\left[U\left(X_T(x) + \frac{X_T(x)}{x}\left(\Delta x(1 + M_T^0) + \sum_{i=1}^{m} q_i(g_i(x) + M_T^i)\right)\right)\right].$$

Applying Lemma 1 to the right-hand side of this inequality, we deduce

$$u(x + \Delta x, q) \geq u(x) + y\Delta x + \langle r(y), q \rangle$$
$$+ \tfrac{1}{2}(\,\Delta x \quad q\,)\, G\left(\begin{matrix}\Delta x \\ q\end{matrix}\right) + o((\Delta x)^2 + \|q\|^2),$$

where the matrix $G$ is given by

$$G_{ij} = -\frac{u'(x)}{x}\mathbb{E}_{\mathbb{R}(x)}[A(X_T(x))(M_T^i + g_i(x))(M_T^j + g_j(x))], \qquad 0 \leq i, j \leq m.$$

Using Lemma 2, we can find, for any $\varepsilon > 0$, the processes $(M^i)_{0 \le i \le m}$ as above such that

$$\|G - G(x)\| \le \varepsilon.$$

This completes the proof.  □

Our next goal is to establish the *upper* bound for $v(y, r)$,

$$v(y + \Delta y, r(y) + \Delta r) \le v(y) - x\Delta y + \tfrac{1}{2} \begin{pmatrix} \Delta y & \Delta r \end{pmatrix} H(y) \begin{pmatrix} \Delta y \\ \Delta r \end{pmatrix}$$

(68)
$$+ o((\Delta y)^2 + \|\Delta r\|^2),$$

where the matrix $H(y)$ is defined in (57) and where we have used the identity $x = -v'(y)$. This will be accomplished in Lemma 6 below and will rely on several auxiliary results.

Denote

(69) $$\mathcal{N}_i^\infty(y) \triangleq \{N \in \mathcal{N}_i^2(y) : N \text{ is bounded}\}, \qquad 0 \le i \le m.$$

LEMMA 4. *Assume that the conditions of Theorem 2 are satisfied. Then for any $0 \le i \le m$, the set $\mathcal{N}_i^\infty(y)$ is nonempty and its closure in $\mathbf{H}^2(\mathbb{R}(x))$ coincides with $\mathcal{N}_i^2(y)$.*

Before proceeding to the proof we note that Assumption 1 is clearly necessary for the sets $\mathcal{N}_i^2(y)$ to be nonempty for all $0 \le i \le m$.

PROOF OF LEMMA 4. To simplify notation, we shall identify a martingale $M$ with its terminal element $M_T$. Let $\mathcal{L}^1(x)$ be the set of stochastic integrals with respect to $S^{X(x)}$ which are martingales under $\mathbb{R}(x)$, that is,

$$\mathcal{L}^1(x) \triangleq \left\{ M = \int H \, dS^{X(x)} : M \text{ is an } \mathbb{R}(x)\text{-martingale} \right\}.$$

It is well known that $\mathcal{L}^1(x)$ is a closed subspace of $\mathbf{L}^1(\mathbb{R}(x))$. From Assumptions 1 and 4, we deduce that no nontrivial linear combination of the random variables $(g_i(x))_{0 \le i \le m}$ is an element of $\mathcal{L}^1(x)$. Hence, $g_i(x)$ (for fixed $0 \le i \le m$) does not belong to the linear subspace generated by $\mathcal{L}^1(x)$ and the random variables $(g_j(x))_{0 \le j \le m, j \ne i}$ (which is closed). We can then apply the Hahn–Banach theorem to strictly separate $g_i(x)$ from this linear subspace. It is easy to see that this provides an element $N_i \in \mathcal{N}_i^\infty(y)$.

Observe now that the *affine* sets $\mathcal{N}_i^2(y)$ are the translations by any element $N^i \in \mathcal{N}_i^\infty(y)$ of the *linear* subspace

$$\mathcal{N}^2(y, g) \triangleq \{N \in \mathcal{N}^2(y) : \mathbb{E}_{\mathbb{R}(x)}[N_T g_j(x)] = 0, 0 \le j \le m\}.$$

Hence, to finish the proof, we need to show that the subspace

$$\mathcal{N}^\infty(y,g) = \{N \in \mathcal{N}^2(y,g) : N \text{ is bounded}\}$$

is dense in $\mathcal{N}^2(y,g)$ [in the $\mathbf{L}^2(\mathbb{R}(x))$ topology].

Denote by $\mathcal{M}^2(x,g)$ the (closed) linear subspace of $\mathbf{L}^2(\mathbb{R}(x))$ generated by $\mathcal{M}^2(x)$ and $(g_i(x))_{0 \le i \le m}$. As $\mathcal{N}^2(y,g)$ is the *complementary* linear subspace of $\mathcal{M}^2(x,g)$, it is enough to show that any element $M \in \mathbf{L}^2(\mathbb{R}(x))$ which is *orthogonal* to $\mathcal{N}^\infty(y,g)$ belongs to $\mathcal{M}^2(x,g)$. Again using the Hahn–Banach theorem, we prove that such an $M$ belongs to the (closed) linear subspace of $\mathbf{L}^1(\mathbb{R}(x))$ generated by $\mathcal{L}^1(x)$ and $(g_i(x))_{0 \le i \le m}$. Since $M$ is square integrable, it turns out that $M \in \mathcal{M}^2(x,g)$.  □

LEMMA 5. *Assume that the conditions of Theorem 2 are satisfied and let $N^i \in \mathcal{N}_i^\infty(y)$ for $0 \le i \le m$. Then there exists $\varepsilon = \varepsilon(N^0, \dots, N^m) > 0$ such that*

$$(70) \qquad \frac{Y(y)}{y}\left(y + \Delta y N^0 + \sum_{i=1}^m \Delta r_i N^i\right) \in \mathcal{Y}(y + \Delta y, r(y) + \Delta r)$$

*for any $(\Delta y, \Delta r) \in \mathbf{R}^{m+1}$ such that $|\Delta y| + \|\Delta r\| \le \varepsilon$.*

PROOF. Since the processes $N^i$ are bounded, there exists $\varepsilon > 0$ such that the process defined in (70) is positive whenever $|\Delta y| + \|\Delta r\| \le \varepsilon$. Fix such $\Delta y$ and $\Delta r$ and denote

$$L \triangleq y + \Delta y N^0 + \sum_{i=1}^m \Delta r_i N^i.$$

Under $\mathbb{R}(x)$, the process $L$ is a positive, bounded martingale and $KL$ is a martingale for any $K \in \mathcal{M}^2(x)$. Using the fact that $S^{X(x)}$ is *sigma-bounded*, we conclude that the product of $L$ on any *nonnegative* stochastic integral with respect to $S^{X(x)}$ is a local martingale (hence, a supermartingale) under $\mathbb{R}(x)$. In particular, as $L_0 = y + \Delta y$, we have

$$\frac{Y(y)}{y} L \in \mathcal{Y}(y + \Delta y).$$

Now let $X \in \mathcal{X}(\tilde{x}, \tilde{q})$ for some $(\tilde{x}, \tilde{q}) \in \mathcal{K}$. Denote

$$M \triangleq \frac{X}{X(x)} x.$$

We claim that $LM$ is a supermartingale under $\mathbb{R}(x)$ and, therefore,

$$(71) \qquad \mathbb{E}_{\mathbb{R}(x)}[L_T M_T] \le \tilde{x}(y + \Delta y).$$

Indeed, from the definition of $\mathcal{X}(\tilde{x}, \tilde{q})$, we deduce that

$$M_T + \sum_{i=1}^{m} \tilde{q}_i g_i(x) \geq 0$$

and from Assumption 4, we deduce the existence of a constant $c > 0$ and a process $N \in \mathcal{M}^2(x)$ such that

$$\sum_{i=1}^{m} \tilde{q}_i |g_i(x)| \leq c + N_T.$$

It follows that $c + N + M$ is a nonnegative stochastic integral with respect to $S^{X(x)}$ and, hence, $L(c + N + M)$ is a supermartingale under $\mathbb{R}(x)$. The supermartingale property for $LM$ now follows from the fact that $L(c + N)$ is a martingale under $\mathbb{R}(x)$.

From the definitions of the sets $\mathcal{N}_i^{\infty}(y)$ and the vector $r(y)$, we deduce that

$$\mathbb{E}_{\mathbb{R}(x)}[L_T g_i(x)] = r_i(y) + \Delta r_i, \qquad 1 \leq i \leq m.$$

This identity and (71) imply that

$$\mathbb{E}\left[\frac{Y_T(y)}{y} L_T \left(X_T + \sum_{i=1}^{m} \tilde{q}_i f_i\right)\right] = \mathbb{E}_{\mathbb{R}(x)}\left[L_T \left(M_T + \sum_{i=1}^{m} \tilde{q}_i g_i(x)\right)\right]$$
$$\leq \tilde{x}(y + \Delta y) + \langle \tilde{q}, r(y) + \Delta r \rangle,$$

thus completing the proof. $\square$

LEMMA 6. *Assume that the conditions of Theorem 2 are satisfied. Then* (68) *holds true.*

PROOF. We closely follow the proof of Lemma 3. Let $N^i \in \mathcal{N}_i^{\infty}(y)$, $0 \leq i \leq m$. By Lemma 5,

$$\frac{Y(y)}{y}\left(y + \Delta y N^0 + \sum_{i=1}^{m} \Delta r_i N^i\right) \in \mathcal{Y}(y + \Delta y, r(y) + \Delta r)$$

for $|\Delta y| + \|\Delta r\|$ sufficiently small, so

$$v(y + \Delta y, r(y) + \Delta r)$$
$$\leq \mathbb{E}\left[V\left(Y_T(y) + \frac{Y_T(y)}{y}\left(\Delta y N_T^0 + \sum_{i=1}^{m} \Delta r_i N_T^i\right)\right)\right].$$

Note now that if $U$ satisfies Assumption 2, then $V$ does also. Hence, we can apply Lemma 1 (or even Lemma 4 from [14]) to the right-hand side of this inequality. This leads to

$$v(y + \Delta y, r(y) + \Delta r) \leq v(y) - x\Delta y + \frac{1}{2} \begin{pmatrix} \Delta y & \Delta r \end{pmatrix} H \begin{pmatrix} \Delta y \\ \Delta r \end{pmatrix}$$
$$+ o((\Delta y)^2 + \|\Delta r\|^2),$$

where the matrix $H$ is given by

$$H_{ij} = \frac{x}{y} \mathbb{E}_{\mathbb{R}(x)}[B(Y_T(y))N_T^i N_T^j], \qquad 0 \leq i, j \leq m.$$

The result now follows from Lemma 4, which allows us to approximate the matrix $H(y)$ by such a matrix $H$.  □

As the functions $u(x, q)$ and $v(y, r)$ satisfy the conjugacy relations (49), their second-order expansions (29) and (52) at the conjugate points $(x, 0)$ and $(y, r(y)) = \partial u(x, 0)$ follow from the bounds (67) and (68) provided that the matrices $-G(x)$ and $H(y)$ are inverse to each other. The proof of this fact relies on the following lemma:

LEMMA 7. *Assume that the conditions of Theorem 2 hold true. Denote by $\alpha(x)$ and $\beta(y)$, respectively, the $(m+1) \times (m+1)$ matrices*

$$(72) \qquad \alpha_{ij}(x) \triangleq \mathbb{E}_{\mathbb{R}(x)}[A(X_T(x))(g_i(x) + M_T^i(x))(g_j(x) + M_T^j(x))],$$

$$(73) \qquad \beta_{ij}(y) \triangleq \mathbb{E}_{\mathbb{R}(x)}[B(Y_T(y))N_T^i(y)N_T^j(y)],$$

*where $0 \leq i, j \leq m$ and where $M(x) = (M^i(x))_{0 \leq i \leq m}$ and $N(y) = (N^i(y))_{0 \leq i \leq m}$ are the solutions to (43) and (54). Then $\alpha(x)$ is positive definite, its inverse equals $\beta(y)$ and (in vector notation)*

$$(74) \qquad A(X_T(x))(g(x) + M_T(x)) = \alpha(x)N_T(y).$$

PROOF. Assumption 1 easily implies that the matrix $\alpha(x)$ has full rank. Using standard arguments from the calculus of variations, we deduce that the solutions $M^i(x)$ to (43) are uniquely determined by the representation

$$A(X_T(x))(g_i(x) + M_T^i(x)) \in a_i + \mathcal{N}^2(y),$$

where $a_i$ is a real number and where we identify the elements of $\mathcal{N}^2(y)$ with their terminal values. Hence, there exists $N = (N^i)_{0 \leq i \leq m}$ such that $N^i - N_0^i \in \mathcal{N}^2(y)$ and

$$(75) \qquad A(X_T(x))(g(x) + M_T(x)) = \alpha(x)N_T.$$

From the definition of $\alpha(x)$, we deduce that

$$\mathbb{E}_{\mathbb{R}(x)}[(g^i(x) + M_T^i(x))N_T^j] = \delta_{ij}$$

and, therefore, $N^i \in \mathcal{N}_i^2(y)$ for $0 \le i \le m$.

Note now that (75) and the fact that $B(Y_T(y)) = 1/A(X_T(x))$ together imply that $B(Y_T(y))N_T^i$ is a linear combination of the elements of the vector $g(x)$ and of the space $\mathcal{M}^2(x)$, which is precisely the characterization of the unique solution to (54). Hence, $N = N(y)$ and (74) follows.

Finally, if we rewrite (74) as

$$g(x) + M_T(x) = \alpha(x)B(Y_T(y))N_T(y)$$

and recall the definition of the matrix $\beta(y)$, then we have that

$$(\alpha(x)\beta(y))_{ij} = \mathbb{E}_{\mathbb{R}(x)}[(g^i(x) + M_T^i(x))N_T^j(y)] = \delta_{ij},$$

where the last equality has already been proven. Hence, the matrices $\alpha(x)$ and $\beta(y)$ are inverse to each other. $\quad\square$

After these preparations, we can proceed to the proofs of Theorems 1 and 2.

PROOFS OF THEOREMS 1 AND 2. First, we note that without any loss of generality, we can assume for the proof of Theorem 1 that Assumption 1 holds true. As before, we denote $y = u'(x)$.

From Lemma 7, we deduce that the matrices $-G(x)$ and $H(y)$ are positive definite and inverse to each other. The same lemma implies the relation (59) between $Z(x)$ and $W(y)$. Further, the second-order approximations (29) [for $u(x,q)$ at $(x,0)$] and (52) [for $v(y,r)$ at $(y,r(y))$] follow from the lower bound (67) for $u(x,q)$ proved in Lemma 3, the upper bound (68) for $v(y,r)$ established in Lemma 6, the conjugacy relations (49) between $u(x,q)$ and $v(y,r)$ and the aforementioned fact that $-G(x)$ is the inverse of $H(y)$.

Note that as the functions $u(x,q)$ and $-v(y,r)$ are concave, the second-order expansions (29) and (52) imply that $G(x)$ and $H(y)$ are the respective matrices of the second derivatives of $u(x,q)$ at $(x,0)$ and $v(y,r)$ at $(y,r(y))$. For $v(y,r)$, the second derivative is defined in the classical sense, as this function is continuously differentiable on its domain $\mathcal{L}$, while for $u(x,q)$ it is understood in the sense of (31).

To complete the proof, we must verify the assertions of the first items of the theorems. By [8], Theorem 2, the solutions $X(x,q)$ to (17) and $Y(y,r)$ to (48) are related by

$$Y_T(y,r) = U'(X_T(x,q) + \langle q, f \rangle), \qquad (y,r) \in \partial u(x,q).$$

From these representations and previously established results, we deduce that the convergences (46) and (56) are, in fact, equivalent. Hereafter, we choose to prove (56).

We proceed along the lines of the proof of Theorem 2 in [14]. Recall that $Y(y, r(y)) = Y(y)$. Let $(y^n, r^n)_{n \geq 1}$ be a sequence in $\mathcal{L}$ that converges to $(y, r(y))$. It is easy to show (see Lemma 3.6 in [13] for a similar result) that

$$(76) \qquad \lim_{n \to \infty} Y_T(y^n, r^n) = Y_T(y, r(y)) = Y_T(y)$$

in probability. Passing, if necessary, to a subsequence, we can assume that this convergence takes place almost surely. In this case,

$$\zeta \triangleq \frac{1}{2} \min \left( Y_T(y), \inf_{n \geq 1} Y_T(y^n, r^n) \right) > 0,$$

$$\eta \triangleq 2 \max \left( Y_T(y), \sup_{n \geq 1} Y_T(y^n, r^n) \right) < \infty$$

and, therefore,

$$(77) \qquad \inf_{\zeta \leq t \leq \eta} V''(t) > 0,$$

where all inequalities hold $\mathbb{P}$-a.s. We define the random variable $\theta$ by

$$\theta \triangleq \min \left( \frac{y X_T(x)}{x Y_T(y)}, \frac{1}{2} \inf_{\zeta \leq t \leq \eta} V''(t) \right).$$

The reasons for this construction of $\theta$ will become clear from the arguments below. Note that, by (77), $\theta$ is strictly positive.

Now let $N^i \in \mathcal{N}_i^\infty(y)$, $0 \leq i \leq m$, where the sets $\mathcal{N}_i^\infty(y)$ were introduced in (69). Define the processes $W = (W^i)_{0 \leq i \leq m}$ and the matrix $H = (H_{ij})_{0 \leq i, j \leq m}$ by

$$W^i \triangleq \frac{Y(y)}{y} N^i,$$

$$H_{ij} \triangleq \frac{x}{y} \mathbb{E}_{\mathbb{R}(x)}[B(Y_T(y)) N_T^i N_T^j].$$

From the construction of the random variable $\theta$, we deduce that

$$(78) \qquad \begin{aligned} \mathbb{E}[\theta \| W_T - W_T(y) \|^2] &\leq \mathbb{E} \left[ \frac{y X_T(x)}{x Y_T(y)} \| W_T - W_T(y) \|^2 \right] \\ &= \mathbb{E}_{\mathbb{R}(x)}[\| N_T - N_T(y) \|^2]. \end{aligned}$$

For $n \geq 1$, denote

$$\begin{aligned} Y^n &\triangleq \frac{Y(y)}{y} \left( y + N^0 \Delta y^n + \sum_{i=1}^m N^i \Delta r_i^n \right) \\ &= Y(y) + W^0 \Delta y^n + \sum_{i=1}^m W^i \Delta r_i^n, \end{aligned}$$

where $\Delta y^n \triangleq y^n - y$ and $\Delta r^n \triangleq r^n - r(y)$. As $(y^n, r^n)$ converge to $(y, r)$ and $N^i$ are bounded, we can assume without loss of generality that for $n \geq 1$,

(79) $$\frac{1}{2} Y(y) \leq Y^n \leq 2Y(y)$$

and (by Lemma 5) that

(80) $$Y^n \in \mathcal{Y}(y^n, r^n).$$

Taking into account (79) and the definition of $\theta$, we deduce that

$$V(Y_T^n) - V(Y_T(y^n, r^n)) \geq V'(Y_T(y^n, r^n))(Y_T^n - Y_T(y^n, r^n))$$
$$+ \theta(Y_T^n - Y_T(y^n, r^n))^2.$$

By Theorem 2 in [8], the terminal value of the optimal solution to (17) at $(x^n, q^n) \triangleq -\partial v(y^n, r^n)$ is given by

$$X_T(x^n, q^n) + \langle q^n, f \rangle = -V'(Y_T(y^n, r^n))$$

and

$$\mathbb{E}[(X_T(x^n, q^n) + \langle q^n, f \rangle) Y_T(y^n, r^n)] = x^n y^n + \langle q^n, r^n \rangle.$$

Together with (80), this implies that

$$\mathbb{E}[V'(Y_T(y^n, r^n))(Y_T^n - Y_T(y^n, r^n))] \geq 0$$

and, therefore, that

(81) $$\mathbb{E}[\theta(Y_T^n - Y_T(y^n, r^n))^2] \leq \mathbb{E}[V(Y_T^n)] - v(y^n, r^n).$$

From Lemma 1 (or Lemma 4 in [14]), we deduce that

$$\mathbb{E}[V(Y_T^n)] = \mathbb{E}\left[ V\left( \frac{Y_T(y)}{y} \left( y + N_T^0 \Delta y^n + \sum_{i=1}^m N_T^i \Delta r_i^n \right) \right) \right]$$
$$= v(y) - x \Delta y^n + \frac{1}{2} \begin{pmatrix} \Delta y^n & \Delta r^n \end{pmatrix} H \begin{pmatrix} \Delta y^n \\ \Delta r^n \end{pmatrix}$$
$$+ o((\Delta y^n)^2 + \|\Delta r^n\|^2).$$

Together with that previously established second-order expansion for $v(y, r)$ at $(y, r(y))$, this implies that

$$\limsup_{n \to \infty} \frac{1}{(\Delta y^n)^2 + \|\Delta r^n\|^2} (\mathbb{E}[V(Y_T^n)] - v(y^n, r^n)) \leq \frac{1}{2} \|H - H(y)\|.$$

Combining the last inequality with (78) and (81), we deduce that

$$\limsup_{n \to \infty} \mathbb{E}\left[ \theta \frac{(Y_T(y^n, r^n) - Y_T(y) - W_T^0(y)\Delta y^n - \sum_{i=1}^m W_T^i(y)\Delta r_i^n)^2}{|\Delta y^n|^2 + \|\Delta r^n\|^2} \right]$$
$$\leq 2 \limsup_{n \to \infty} \mathbb{E}\left[ \theta \left( \frac{(Y_T(y^n, r^n) - Y_T^n)^2}{|\Delta y^n|^2 + \|\Delta r^n\|^2} + \|W_T(y) - W_T\|^2 \right) \right]$$
$$\leq \|H - H(y)\| + 2\mathbb{E}_{\mathbb{R}(x)}[\|N_T - N_T(y)\|^2].$$

Lemma 4 allows us to make the right-hand side an *arbitrary* small number. The result, that is, the convergence (56), now follows from the fact that $\theta > 0$. □

## 4. The computation of the sensitivity parameters $p'(x)$ and $D(x)$.

As before, we fix $x > 0$ and denote $y = u'(x)$. We assume throughout this section that the conditions of Theorem 1 hold true. The sensitivity parameters $p'(x)$ and $D(x)$ will be computed in terms of the first- and second-order derivatives of the value function $u(x, q)$ at $(x, 0)$. To facilitate the subsequent "qualitative" study of $p'(x)$ and $D(x)$, we start by giving more convenient expressions for those elements of the matrix $G(x)$ [defined in (47)] that contain index 0 and are thus related to the differentiability with respect to the liquid wealth $x$. First, clearly,

$$G_{00}(x) = u''(x),$$

where $u''(x)$ is the second derivative of the value function $u$ defined in (18).

To perform further computations, we introduce the derivative process $Y'(y)$ of the solution $Y(y)$ to (22). According to [14], Theorem 1, such a process $Y'(y)$ is uniquely defined by the terminal condition

$$Y'_T(y) = \lim_{\Delta y \to 0} \frac{Y_T(y + \Delta y) - Y_T(y)}{\Delta y},$$

where the convergence takes place in $\mathbb{P}$-probability, and by the martingale property of $X(x)Y'(y)$. Moreover, if we denote

$$L(y) \triangleq \frac{y}{Y(y)} Y'(y),$$

then $L(y)$ belongs to $1 + \mathcal{N}^2(y)$, where $\mathcal{N}^2(y)$ is defined in (53) and solves the following optimization problem:

$$c(y) \triangleq \min_{L \in 1 + \mathcal{N}^2(y)} \mathbb{E}_{\mathbb{R}(x)}[B(Y_T(y))(L_T)^2].$$

From the duality relation (21), we deduce that

$$U''(X_T(x))X'_T(x) = u''(x)Y'_T(y).$$

It follows that

$$\begin{aligned}
G_{0i}(x) = G_{i0}(x) &= \mathbb{E}[U''(X_T(x))X'_T(x)(Z^i_T(x) + f_i)] \\
&= u''(x)\mathbb{E}[Y'_T(y)(Z^i_T(x) + f_i)] = u''(x)\mathbb{E}[Y'_T(y)f_i] \\
&= u''(x)\widetilde{p}_i(x), \qquad 1 \le i \le m,
\end{aligned}$$

where we used the fact that $Y'(y)Z^i(x)$ is a $\mathbb{P}$-martingale starting at 0 and the notation

$$(82) \qquad \widetilde{p}_i(x) \triangleq \mathbb{E}[Y'_T(y)f_i], \qquad 1 \le i \le m.$$

THEOREM 3. *Assume that the conditions of Theorem* 1 *are satisfied. Then the marginal utility-based prices satisfy* (28) *and the corresponding sensitivity parameters* $p'(x)$ *and* $D(x)$ *equal*

$$p_i'(x) = \frac{u''(x)}{u'(x)}(\widetilde{p}_i(x) - p_i(x)), \tag{83}$$

$$D_{ij}(x) = \frac{1}{u'(x)}G_{ij}(x) - \frac{u''(x)}{u'(x)}p_i(x)\widetilde{p}_j(x), \tag{84}$$

*where* $p_i(x)$ *and* $\widetilde{p}_i(x)$ *are defined in* (25) *and* (82) *and* $G_{ij}(x)$ *is given by* (47), $1 \leq i, j \leq m$.

PROOF. Taking into account the characterization (24) of the marginal utility-based prices the relations (83) and (84), are the results of direct computations based on Theorem 1 and Lemma 8 below. □

LEMMA 8. *Let* $D \subset \mathbb{R}^N$ *be an open convex set and let* $w : D \to \mathbb{R}$ *be a convex function. Assume that for some* $x_0 \in D$, *the function* $w$ *admits the second-order expansion*

$$w(x) = w(x_0) + \langle p_0, x - x_0 \rangle + \tfrac{1}{2}\langle A(x - x_0), x - x_0 \rangle + o(\|x - x_0\|^2), \tag{85}$$

*for some* $p_0 \in \mathbb{R}^N$ *and a symmetric* $N \times N$ *matrix* $A$. *Assume also that* $f(p)$ *is a continuously differentiable function in a neighborhood of* $p_0$. *We then have the following chain rule for* $f(\partial w)$:

$$\lim_{x \to x_0} \sup_{p \in \partial w(x)} \frac{|f(p) - f(p_0) - \langle \nabla f(p_0), A(x - x_0) \rangle|}{\|x - x_0\|} = 0.$$

PROOF. The proof of the lemma is based on the mean value theorem,

$$f(p) - f(p_0) = \langle \nabla f(\tilde{p}), p - p_0 \rangle,$$

where $\tilde{p} = \lambda p + (1 - \lambda)p_0$ for some $0 \leq \lambda \leq 1$, as on well as the well-known observation that once the second-order expansion (85) has been established then

$$\lim_{x \to x_0} \sup_{p \in \partial w(x)} \frac{\|p - p_0 - A(x - x_0)\|}{\|x - x_0\|} = 0. \tag{86}$$

For the convenience of the reader, we include here a short proof of (86). We can assume without loss of generality that in (85), $x_0 = p_0 = 0$ and $w(0) = 0$. Then we have

$$w(x) = \tfrac{1}{2}\langle Ax, x \rangle + o(\|x\|^2). \tag{87}$$

From the definition of the subgradient $\partial w(x)$ and (87), it follows that for any $p \in \partial w(x)$ and $z \in \mathbf{R}^N$,

$$\langle p - Ax, z \rangle = \langle p, z \rangle - \langle Ax, z \rangle \leq w(x + z) - w(x) - \langle Ax, z \rangle$$
$$\leq \tfrac{1}{2} \langle Az, z \rangle + o(\|x + z\|^2) + o(\|x\|^2).$$

Fix $t > 0$ and choose $z = t\|x\| \frac{p - Ax}{\|p - Ax\|}$. We have

$$t\|x\| \|p - Ax\| \leq \tfrac{1}{2} \|A\| t^2 \|x\|^2 + o(2(t^2 + 1)\|x\|^2) + o(\|x\|^2).$$

Dividing by $t\|x\|^2$ and letting $x \to 0$, we obtain

$$\limsup_{x \to 0} \sup_{p \in \partial w(x)} \frac{\|p - Ax\|}{\|x\|} \leq \frac{1}{2} \|A\| t.$$

As $t$ is any positive number, this implies (86) (for the case $x_0 = p_0 = 0$). The proof of the lemma is now complete. $\quad\square$

## 5. Risk-tolerance wealth processes.

We recall that the quantity

$$-\frac{U'(x)}{U''(x)} = \frac{x}{A(x)}$$

is called the *absolute risk-tolerance coefficient* of the utility function $U$ computed at $x > 0$. To the best of our knowledge, the following concept has not been defined in the literature previously:

DEFINITION 2. Let $x > 0$ and assume that the solution $X(x)$ to (18) exists. We call a process $R(x)$ the *risk-tolerance wealth process* if it is a maximal wealth process and

$$(88) \qquad R_T(x) = -\frac{U'(X_T(x))}{U''(X_T(x))}.$$

In other words, $R(x)$ is the replication process for the random payoff defined in the right-hand side of (88).

The risk-tolerance wealth process $R(x)$ may very well not exist. In fact, as we shall see in Section 6, the existence of $R(x)$ is closely related to the validity of the important qualitative properties (3)–(6) for marginal utility-based prices stated in Section 1. Motivated by this observation, in this section we present different criteria for $R(x)$ to be well defined.

We start with a basic result where the conditions are given in terms of the value functions of the optimization problems (18) and (22).

THEOREM 4. *Assume that conditions* (11) *and* (19) *and Assumption* 2 *hold true. Fix $x > 0$ and denote $y = u'(x)$. Then the following assertions are equivalent:*

1. *The risk-tolerance wealth process $R(x)$ exists.*
2. *The value function $u(x)$ in* (18) *is two-times differentiable at $x$ and*

$$\frac{(u'(x))^2}{u''(x)} = \mathbb{E}\left[\frac{(U'(X_T(x)))^2}{U''(X_T(x))}\right], \tag{89}$$

   *where $X(x)$ is the solution to* (18).
3. *The value function $v(y)$ in* (22) *is two-times differentiable at $y = u'(x)$ and*

$$y^2 v''(y) = \mathbb{E}[Y_T^2(y)V''(Y_T(y))], \tag{90}$$

   *where $Y(y)$ is the solution to* (22).

*In addition, if these assertions hold true, then the initial value of $R(x)$ is given by the absolute risk-tolerance coefficient of $u(x)$,*

$$R_0(x) = -\frac{u'(x)}{u''(x)}, \tag{91}$$

*the product $R(x)Y(y)$ is a uniformly integrable martingale and*

$$\lim_{\varepsilon \to 0} \frac{X_T(x+\varepsilon) - X_T(x)}{\varepsilon} = \frac{R_T(x)}{R_0(x)}, \tag{92}$$

$$\lim_{\varepsilon \to 0} \frac{Y_T(y+\varepsilon) - Y_T(y)}{\varepsilon} = \frac{Y_T(y)}{y}, \tag{93}$$

*where the convergence takes place in $\mathbb{P}$-probability.*

Remark 2. The convergences (92) and (93) imply that

$$\frac{1}{R_0(x)} R(x) = X'(x), \tag{94}$$

$$\frac{1}{y} Y(y) = Y'(y), \tag{95}$$

where $X'(x)$ and $Y'(y)$ are the derivatives of the solutions to (18) and (22) introduced in the previous sections. In particular, we see that if the risk-tolerance wealth process exists, then it represents the behavior of an economic agent with respect to an infinitesimally small increase of the initial wealth.

Note that, formally speaking, the existence of $X'(x)$ and $Y'(y)$, as well as $u''(x)$ and $v''(y)$, requires, in addition to the conditions of Theorem 4, the validity of Assumption 3 (concerning the sigma-boundedness of $S^{X(x)}$). We refer the reader to [14] for the exact statement and a counterexample. In fact, the absence of this assumption represents the main difficulty for the proof of the theorem.

PROOF OF THEOREM 4.   First, the equivalence of (89) and (90) is an easy consequence of the duality relations between the value functions and the solutions to the optimization problems (18) and (22).

Assume now that the risk-tolerance wealth process $R(x)$ exists. Then $R(x)/X(x)$ replicates the random variable $1/A(X_T(x))$ in the model where $X(x)$ is chosen as a numéraire. Since $A(X_T(x))$ is bounded away from zero, we can conclude that $R(x)/X(x)$ is a *bounded* martingale under the probability measure $\mathbb{R}(x)$ introduced in (40). It follows that $R(x)Y(y)$ is a martingale under $\mathbb{P}$.

For sufficiently small $\Delta x$ and $\Delta y$, we deduce from Lemma 4 in [14] that

$$u(x + \Delta x) \geq \mathbb{E}\left[U\left(X_T(x) + \Delta x \frac{R_T(x)}{R_0(x)}\right)\right]$$

$$= u(x) + y\Delta x + \frac{1}{2}(\Delta x)^2 \mathbb{E}\left[U''(X_T(x))\left(\frac{R_T(x)}{R_0(x)}\right)^2\right] + o((\Delta x)^2),$$

$$v(y + \Delta y) \leq \mathbb{E}\left[V\left(Y_T(y) + \Delta y \frac{Y_T(y)}{y}\right)\right]$$

$$= v(y) - x\Delta y + \frac{1}{2}(\Delta y)^2 \mathbb{E}\left[V''(Y_T(y))\left(\frac{Y_T(y)}{y}\right)^2\right] + o((\Delta y)^2).$$

Since

$$U''(X_T(x))R_T(x) = -Y_T(y),$$

$$V''(Y_T(y))Y_T(y) = R_T(x)$$

and $R(x)Y(y)$ is a uniformly integrable martingale, we conclude that

$$\mathbb{E}\left[U''(X_T(x))\left(\frac{R_T(x)}{R_0(x)}\right)^2\right]\mathbb{E}\left[V''(Y_T(y))\left(\frac{Y_T(y)}{y}\right)^2\right] = -1.$$

It follows that the second derivatives $u''(x)$ and $v''(y)$ are well defined and that

$$u''(x) = \mathbb{E}\left[U''(X_T(x))\left(\frac{R_T(x)}{R_0(x)}\right)^2\right],$$

$$v''(y) = \mathbb{E}\left[V''(Y_T(y))\left(\frac{Y_T(y)}{y}\right)^2\right].$$

It is easy to see that the last expression coincides with the one given in (90).

Assume now that (90) holds true. We denote by $\mathcal{N}^\infty(y)$ the family of semimartingales $N$ such that for some $\varepsilon = \varepsilon(N) > 0$,

$$Y(y)(1 + \varepsilon N) \in \mathcal{Y}(y) \quad \text{and} \quad Y(y)(1 - \varepsilon N) \in \mathcal{Y}(y).$$

For $N \in \mathcal{N}^\infty(y)$, we have that for sufficiently small $|\Delta y|$,

$$Y(y)\left(1 + \frac{\Delta y}{y}(1 + N)\right) \in \mathcal{Y}(y + \Delta y)$$

and, hence,

$$v(y + \Delta y) \leq \mathbb{E}\left[V\left(Y_T(y)\left(1 + \frac{\Delta y}{y}(1 + N_T)\right)\right)\right].$$

From (90) and [14], Lemma 4, we deduce that

$$y^2 v''(y) = \mathbb{E}[V''(Y_T(y))(Y_T(y))^2] \leq \mathbb{E}[V''(Y_T(y))(Y_T(y)(1 + N_T))^2].$$

As $N$ was an arbitrary element of the linear space $\mathcal{N}^\infty(y)$, we deduce that

$$\mathbb{E}[V''(Y_T(y))(Y_T(y))^2 N_T] = 0, \qquad N \in \mathcal{N}^\infty(y),$$

or, equivalently, using the probability measure $\mathbb{R}(x)$ introduced in (40) and the function $B$ defined in (35),

$$\mathbb{E}_{\mathbb{R}(x)}[B(Y_T(y))N_T] = 0, \qquad N \in \mathcal{N}^\infty(y).$$

The above identity can also be written as

$$(96) \qquad \mathbb{E}_{\mathbb{Q}}[B(Y_T(y))] = \mathbb{E}_{\mathbb{R}(x)}[B(Y_T(y))], \quad \mathbb{Q} \in \widetilde{\mathcal{R}}^0,$$

where

$$\widetilde{\mathcal{R}}^0 \triangleq \left\{\mathbb{Q} \sim \mathbb{P} : \frac{d\mathbb{Q}}{d\mathbb{R}(x)} = 1 + N_T, N \in \mathcal{N}^\infty(y)\right\}.$$

By Theorem B.1 given in Appendix B, we deduce from (96) that the bounded random variable

$$B(Y_T(y)) = \frac{1}{A(X_T(x))}$$

is replicable in the financial model where $X(x)$ is chosen as the numéraire. Hence, $X_T(x)/A(X_T(x))$ is replicable under the original numéraire and the risk-tolerance wealth process exists.

To conclude the proof of the theorem, we need to verify assertions (91), (92) and (93). Expression (91) for the initial value of $R(x)$ is an immediate consequence of (89) and the aforementioned fact that $R(x)Y(y)$ is a uniformly integrable martingale. Further, if (93) holds true, then (92) easily follows:

$$\lim_{\varepsilon \to 0} \frac{X_T(x + \varepsilon) - X_T(x)}{\varepsilon} = \lim_{\varepsilon \to 0} \frac{-V'(Y_T(u'(x + \varepsilon)) + V'(Y_T(u'(x)))}{\varepsilon}$$

$$= -V''(Y_T(y))\frac{Y_T(y)}{y}u''(x) = \frac{R_T(x)}{R_0(x)}.$$

To prove (93), we proceed along the lines of the proof of [14], Theorem 2. Note that similar arguments were used in the proofs of the first items of Theorems 1 and 2. Let $(\varepsilon_n)_{n \geq 1}$ be a sequence of real numbers converging to 0. By Lemma 3.6 in [13],

$$\lim_{n \to \infty} Y_T(y + \varepsilon_n) = Y_T(y)$$

in probability. Passing, if necessary, to a subsequence, we can assume that this convergence takes place almost surely. In this case, there is a strictly positive random variable $\theta$ such that $\theta \leq \frac{1}{2} V''(t)$ for

$$\frac{1}{2} \min \left\{ Y_T(y), \inf_{n \geq 1} Y_T(y + \varepsilon_n) \right\} \leq t \leq 2 \max \left\{ Y_T(y), \sup_{n \geq 1} Y_T(y + \varepsilon_n) \right\}.$$

Denoting

$$\zeta_n \triangleq Y_T(y) \left( 1 + \frac{\varepsilon_n}{y} \right) - Y_T(y + \varepsilon_n),$$

we deduce that for $\frac{1}{2} \leq 1 + \frac{\varepsilon_n}{y} \leq 2$,

$$V \left( Y_T(y) \left( 1 + \frac{\varepsilon_n}{y} \right) \right) - V(Y_T(y + \varepsilon_n)) \geq V'(Y_T(y + \varepsilon_n)) \zeta_n + \theta(\zeta_n)^2.$$

As $-V'(Y_T(y + \varepsilon_n)) = X_T(-v'(y + \varepsilon_n))$, we obtain that

$$\mathbb{E}[V'(Y_T(y + \varepsilon_n)) \zeta_n] \geq 0$$

and, hence, that

$$\mathbb{E}[\theta(\zeta_n)^2] \leq \mathbb{E} \left[ V \left( Y_T(y) \left( 1 + \frac{\varepsilon_n}{y} \right) \right) \right] - v(y + \varepsilon_n).$$

Taking into account (90), we deduce from Lemma 4 in [14] that

$$\mathbb{E} \left[ V \left( Y_T(y) \left( 1 + \frac{\varepsilon_n}{y} \right) \right) \right] = v(y) + v'(y) \varepsilon_n + \frac{1}{2} v''(y) \varepsilon_n^2 + o(\varepsilon_n^2).$$

It follows that

$$\lim_{n \to \infty} \frac{1}{\varepsilon_n^2} \mathbb{E}[\theta(\zeta_n)^2] = 0,$$

thus proving (93). The proof is thus complete.   $\square$

The "global" version of Theorem 4 that requires the existence of $R(x)$ for all $x > 0$ has a particularly elegant formulation in terms of the solutions $Y(y)$ to the dual problem (22).

THEOREM 5.  *Assume that conditions* (11) *and* (19) *and Assumption* 2 *hold true. Then the following assertions are equivalent:*

1. *The risk-tolerance wealth process $R(x)$ exists for any $x > 0$.*
2. *The optimal solution $Y(y)$ to (22) is uniquely defined up to a multiplicative constant, that is,*

$$(97) \qquad Y(y) = yY(1), \qquad y > 0.$$

PROOF.   $2 \Longrightarrow 1$; By (97), we have

$$v(y) = \mathbb{E}[V(yY_T(1))].$$

From [14], Lemma 4, we deduce that

$$v''(y) = \mathbb{E}[V''(yY_T(1))(Y_T(1))^2] = \mathbb{E}\left[V''(Y_T(y))\left(\frac{Y_T(y)}{y}\right)^2\right].$$

This means that for every $y > 0$, the relation (90) holds true. According to Theorem 4, the risk-tolerance wealth process then exists for any $x > 0$.

$1 \Longrightarrow 2$: The proof is easy for the case of a finite probability space. Indeed, from Theorem 4 [see (93)], we deduce that

$$(98) \qquad \frac{d}{dy}\left(\frac{Y_T(y)}{y}\right) = 0, \qquad y > 0,$$

where the derivative is considered in the sense of convergence in probability. If $|\Omega| < \infty$, this then implies (97).

If $\Omega$ is infinite, then (98) is no longer sufficient for (97). In this case, the proof is more involved and is divided into a series of lemmas. For the convenience of future references, we start with the following, easy corollary of Assumption 2.

LEMMA 9.   *Assume that the utility function $U = U(x)$ satisfies (33). Then for any constant $0 < a < 1$, there exist constants $0 < b_1 < b_2 < \infty$ such that*

$$b_1|U''(x)| \leq \min_{|z-x| \leq ax} |U''(z)| \leq \max_{|z-x| \leq ax} |U''(z)| \leq b_2|U''(x)|, \qquad x > 0.$$

PROOF.   The result is an immediate consequence of (33) and Lemma 3 in [14].   $\square$

Of course, similar inequalities also hold true for the conjugate function $V = V(y)$. Below, we shall use this fact without further comment.

LEMMA 10.   *Assume that the conditions of Theorem 5 hold true and that the risk-tolerance wealth process $R(x)$ exists for any $x > 0$. Let $0 < z_1 < z_2 < \infty$ and denote for $n \geq 1$ and $0 \leq i \leq n$,*

$$(99) \qquad y_{i,n} = z_1 + i\frac{z_2 - z_1}{n}, \qquad \varepsilon_{i,n} = \frac{z_2 - z_1}{ny_{i,n}}.$$

*Then*

$$\lim_{n \to \infty} n \sum_{i=0}^{n-1} \left[ \mathbb{E}[V(Y_T(y_{i,n})(1 + \varepsilon_{i,n}))] - v(y_{i+1,n}) \right] = 0. \tag{100}$$

PROOF. As

$$\lim_{n \to \infty} \frac{n}{z_2 - z_1} \sum_{i=1}^{n-1} [v(y_{i-1,n}) - 2v(y_{i,n}) + v(y_{i+1,n})]$$

$$= \lim_{n \to \infty} \frac{n}{z_2 - z_1} \left( \left( v(z_2) - v\left( z_2 - \frac{z_2 - z_1}{n} \right) \right) \right.$$

$$\left. - \left( v\left( z_1 + \frac{z_2 - z_1}{n} \right) - v(z_1) \right) \right)$$

$$= v'(z_2) - v'(z_1),$$

we deduce that the convergence (100) holds true if

$$\lim_{n \to \infty} \frac{n}{z_2 - z_1} \sum_{i=1}^{n-1} \mathbb{E}[V(Y_T(y_{i,n})(1 - \varepsilon_{i,n})) - 2V(Y_T(y_{i,n}))$$

$$+ V(Y_T(y_{i,n})(1 + \varepsilon_{i,n}))] \tag{101}$$

$$= v'(z_2) - v'(z_1).$$

For $0 < \varepsilon < 1$, define the following functions on $(0, \infty)$:

$$F_\varepsilon(y) \triangleq y^2 \min_{|z-y| \le \varepsilon y} V''(z), \qquad G_\varepsilon(y) \triangleq y^2 \max_{|z-y| \le \varepsilon y} V''(z).$$

From Lemma 9, we deduce the existence of $c > 0$ such that

$$0 \le F_\varepsilon(y) \le y^2 V''(y) \le G_\varepsilon(y) \le c y^2 V''(y), \qquad y > 0. \tag{102}$$

It is easy to see that

$$F_\varepsilon(y) \le \frac{1}{\varepsilon^2} (V(y(1 - \varepsilon)) - 2V(y) + V(y(1 + \varepsilon))) \tag{103}$$

$$\le G_\varepsilon(y), \qquad y > 0.$$

By Theorem 4, we have

$$y^2 v''(y) = \mathbb{E}[Y_T^2(y) V''(Y_T(y))], \qquad y > 0.$$

Taking into account Assumption 2 and the uniform integrability of the family $\{X_T(-v'(y))Y_T(y)\}_{z_1 \le y \le z_2}$, we obtain that $v''$ is continuous and, hence,

*Riemann integrable.* It follows that

$$v'(z_2) - v'(z_1) = \int_{z_1}^{z_2} v''(y)\,dy = \lim_{n\to\infty} \frac{z_2 - z_1}{n} \sum_{i=1}^{n} v''(y_{i,n})$$

$$= \lim_{n\to\infty} \frac{z_2 - z_1}{n} \sum_{i=1}^{n} \mathbb{E}\left[V''(Y_T(y_{i,n}))\left(\frac{Y_T(y_{i,n})}{y_{i,n}}\right)^2\right].$$

Using (102), (103) and the result immediately above, we deduce that the convergence (101) holds true if

$$(104) \qquad \lim_{n\to\infty} \frac{1}{n} \sum_{i=1}^{n} \mathbb{E}\left[(G_{\varepsilon_{i,n}}(Y_T(y_{i,n})) - F_{\varepsilon_{i,n}}(Y_T(y_{i,n})))\frac{1}{y_{i,n}^2}\right] = 0.$$

For $z_1 \leq y \leq z_2$, denote

$$Q_n(y) \triangleq \sum_{i=1}^{n}\left[(G_{\varepsilon_{i,n}}(Y_T(y_{i,n})) - F_{\varepsilon_{i,n}}(Y_T(y_{i,n})))\frac{1}{y_{i,n}^2}I_{\{y_{i-1,n}\leq y<y_{i,n}\}}\right].$$

Since

$$\lim_{y_n\to y} Y_T(y_n) = Y_T(y)$$

in probability and $V'' = V''(y)$ is continuous, we have

$$\lim_{n\to\infty} Q_n(y) = 0$$

in probability for any $y \in [z_1, z_2]$. Furthermore, by (102) and Assumption 2 [see (34)], there exists $C > 0$ such that

$$Q_n(y) \leq \frac{c}{z_1^2} V''(Y_T(y_n^*))(Y_T(y_n^*))^2 \leq CX_T(-v'(y_n^*))Y_T(y_n^*)$$

for some $y_n^* \in [z_1, z_2]$. Again taking into account the fact that the family of random variables $\{X_T(-v'(y))Y_T(y)\}_{z_1\leq y\leq z_2}$ is uniformly integrable, we conclude that

$$\lim_{n\to\infty} \mathbb{E}[Q_n(y)] = 0$$

for any $y \in [z_1, z_2]$. We can easily see that

$$\mathbb{E}[Q_n(y)] \leq C \sup_{z_1\leq z\leq z_2} v''(z),$$

so we can use the dominated convergence theorem to obtain

$$\lim_{n\to\infty} \int_{z_1}^{z_2} \mathbb{E}[Q_n(y)]\,dy = 0.$$

This is equivalent to (104), so the proof of the lemma is complete. $\quad\square$

LEMMA 11. *Assume the conditions and the notation of Lemma* 10. *Then*

$$\lim_{n\to\infty}\sum_{i=0}^{n-1}\mathbb{E}\bigg[X_T(x_{i,n})\bigg(Y_T(y_{i,n})\frac{y_{i+1,n}}{y_{i,n}}-Y_T(y_{i+1,n})\bigg)$$
$$\times I_{\{Y_T(y_{i,n})>Y_T(y_{i+1,n})\}}\bigg]=0,\tag{105}$$

*where* $x_{i,n}\triangleq-v'(y_{i,n})$.

PROOF. From Assumption 2 [see (34)] and Lemma 3 in [14], we deduce the existence of a constant $a>0$ such that for any $y>0$,

$$\frac{1}{2}V''(t)\geq a\frac{X_T(-v'(y))}{Y_T(y)},\qquad t<2Y_T(y).$$

It follows that for any $0<y<z<\infty$ such that $z<2y$, we have

$$V\bigg(Y_T(y)\frac{z}{y}\bigg)-V(Y_T(z))$$
$$\geq V'(Y_T(z))\bigg(Y_T(y)\frac{z}{y}-Y_T(z)\bigg)$$
$$+a\frac{X_T(-v'(y))}{Y_T(y)}\bigg(Y_T(y)\frac{z}{y}-Y_T(z)\bigg)^2 I_{\{Y_T(y)>Y_T(z)\}}.$$

As $-V'(Y_T(z))=X_T(-v'(z))$, we obtain that

$$\mathbb{E}\bigg[V'(Y_T(z))\bigg(Y_T(y)\frac{z}{y}-Y_T(z)\bigg)\bigg]\geq 0.$$

It follows that

$$\mathbb{E}\bigg[V\bigg(Y_T(y)\frac{z}{y}\bigg)\bigg]-v(z)$$
$$\geq a\mathbb{E}\bigg[\frac{X_T(-v'(y))}{Y_T(y)}\bigg(Y_T(y)\frac{z}{y}-Y_T(z)\bigg)^2 I_{\{Y_T(y)>Y_T(z)\}}\bigg].$$

From Lemma 10, we deduce that

$$\lim_{n\to\infty}n\sum_{i=0}^{n-1}\mathbb{E}\bigg[\frac{X_T(x_{i,n})}{Y_T(y_{i,n})}\bigg(Y_T(y_{i,n})\frac{y_{i+1,n}}{y_{i,n}}-Y_T(y_{i+1,n})\bigg)^2 I_{\{Y_T(y_{i,n})>Y_T(y_{i+1,n})\}}\bigg]=0.$$

The result now follows from the following chain of inequalities, where $\mathbb{R}(x)$ is the probability measure defined in (40),

$$b\triangleq\min_{z_1\leq y\leq z_2}|yv'(y)|,\qquad c\triangleq\max_{z_1\leq y\leq z_2}|yv'(y)|$$

and where, in the second step, we use the Cauchy inequality:

$$n \sum_{i=0}^{n-1} \mathbb{E}\left[\frac{X_T(x_{i,n})}{Y_T(y_{i,n})}\left(Y_T(y_{i,n})\frac{y_{i+1,n}}{y_{i,n}} - Y_T(y_{i+1,n})\right)^2 I_{\{Y_T(y_{i,n}) > Y_T(y_{i+1,n})\}}\right]$$

$$\geq b \sum_{i=0}^{n-1} \frac{1}{n} \mathbb{E}_{\mathbb{R}(x_{i,n})}\left[n^2\left(\frac{y_{i+1,n}}{y_{i,n}} - \frac{Y_T(y_{i+1,n})}{Y_T(y_{i,n})}\right)^2 I_{\{Y_T(y_{i,n}) > Y_T(y_{i+1,n})\}}\right]$$

$$\geq b \left(\sum_{i=0}^{n-1} \mathbb{E}_{\mathbb{R}(x_{i,n})}\left[\left(\frac{y_{i+1,n}}{y_{i,n}} - \frac{Y_T(y_{i+1,n})}{Y_T(y_{i,n})}\right) I_{\{Y_T(y_{i,n}) > Y_T(y_{i+1,n})\}}\right]\right)^2$$

$$\geq \frac{b}{c^2} \left(\sum_{i=0}^{n-1} \mathbb{E}\left[X_T(x_{i,n})\left(Y_T(y_{i,n})\frac{y_{i+1,n}}{y_{i,n}} - Y_T(y_{i+1,n})\right)\right.\right.$$

$$\left.\left. \times I_{\{Y_T(y_{i,n}) > Y_T(y_{i+1,n})\}}\right]\right)^2. \qquad \square$$

LEMMA 12. *Assume that the conditions of Lemma* 10 *hold true and let* $0 < z_1 < z_2 < \infty$. *Then*

$$Y_T(z_1) \leq Y_T(z_2).$$

PROOF. For $n \geq 1$ and $0 \leq i \leq n$, define $y_{i,n}$ and $\varepsilon_{i,n}$ by (99) and denote $x_{i,n} = -v'(y_{i,n})$. As $V$ is convex and $V'(Y_T(y_{i,n})) = -X_T(x_{i,n})$, we deduce that

$$V(Y_T(z_1)) - V(Y_T(z_2))$$

$$= \sum_{i=0}^{n-1} [V(Y_T(y_{i,n})) - V(Y_T(y_{i+1,n}))]$$

$$\geq -\sum_{i=0}^{n-1} [X_T(x_{i+1,n})(Y_T(y_{i,n}) - Y_T(y_{i+1,n}))]$$

(106)

$$\geq -\sum_{i=0}^{n-1} [(X_T(x_{i+1,n}) - X_T(x_{i,n}))(Y_T(y_{i,n}) - Y_T(y_{i+1,n}))]$$

$$\quad - \sum_{i=0}^{n-1} [X_T(x_{i,n})(Y_T(y_{i,n})(1 + \varepsilon_{i,n}) - Y_T(y_{i+1,n}))$$

$$\quad\quad\quad \times I_{\{Y_T(y_{i,n}) > Y_T(y_{i+1,n})\}}].$$

As

$$\mathbb{E}[(X_T(x_{i+1,n}) - X_T(x_{i,n}))(Y_T(y_{i,n}) - Y_T(y_{i+1,n}))]$$

$$\leq |x_{i+1,n} - x_{i,n}||y_{i+1,n} - y_{i,n}|$$

$$= \frac{z_2 - z_1}{n}(v'(y_{i+1,n}) - v'(y_{i,n})),$$

we deduce that

$$\sum_{i=0}^{n-1} \mathbb{E}[(X_T(x_{i+1,n}) - X_T(x_{i,n}))(Y_T(y_{i,n}) - Y_T(y_{i+1,n}))]$$

$$\leq \frac{z_2 - z_1}{n}(v'(z_2) - v'(z_1))$$

and since

$$(X_T(x_{i+1,n}) - X_T(x_{i,n}))(Y_T(y_{i,n}) - Y_T(y_{i+1,n})) \geq 0,$$

the first term in the right-hand side of (106) converges to zero in probability as $n \to \infty$. By Lemma 11, the second term in the right-hand side of (106) also converges to zero in probability as $n \to \infty$. It follows that

$$V(Y_T(z_1)) - V(Y_T(z_2)) \geq 0.$$

As $V$ is strictly decreasing, this implies the assertion in the statement of the lemma.  $\square$

LEMMA 13.  *Assume that the conditions of Lemma 10 hold true. Let $a > 0$, denote $b \triangleq u'(a)$ and*

$$\eta(y) \triangleq \sqrt{\frac{X_T(a)}{Y_T(b)}}\frac{Y_T(y)}{y}, \qquad y > 0.$$

*Then for $y < b$,*

(107) 
$$\lim_{\varepsilon \to 0}\frac{1}{\varepsilon}\|\eta(y + \varepsilon) - \eta(y)\|_{\mathbf{L}^2} = 0.$$

PROOF.  From Assumption 2 [see (34)] and Lemma 3 in [14], we deduce the existence of a constant $c > 0$ such that

$$\frac{1}{2}V''(t) \geq c\frac{X_T(a)}{Y_T(b)}, \qquad t < 2Y_T(b).$$

By Lemma 12 for $0 < y, y + \varepsilon < b$, $\varepsilon < y$, we have

$$V\left(Y_T(y)\left(1 + \frac{\varepsilon}{y}\right)\right) - V(Y_T(y + \varepsilon))$$

$$\geq V'(Y_T(y + \varepsilon))\left(Y_T(y)\left(1 + \frac{\varepsilon}{y}\right) - Y_T(y + \varepsilon)\right)$$

$$+ c\frac{X_T(a)}{Y_T(b)}\left(Y_T(y)\left(1 + \frac{\varepsilon}{y}\right) - Y_T(y + \varepsilon)\right)^2.$$

As $-V'(Y_T(y + \varepsilon)) = X_T(-v'(y + \varepsilon))$, it follows that

$$\mathbb{E}\Big[V'(Y_T(y + \varepsilon))\Big(Y_T(y)\Big(1 + \frac{\varepsilon}{y}\Big) - Y_T(y + \varepsilon)\Big)\Big] \geq 0.$$

Recalling the notation for the random variables $\eta$, we deduce that

$$\mathbb{E}\Big[V\Big(Y_T(y)\Big(1 + \frac{\varepsilon}{y}\Big)\Big)\Big] - v(y + \varepsilon) \geq c(y + \varepsilon)^2 (\|\eta(y + \varepsilon) - \eta(y)\|_{\mathbf{L}^2})^2.$$

The result now follows from the expansion

$$\mathbb{E}\Big[V\Big(Y_T(y)\Big(1 + \frac{\varepsilon}{y}\Big)\Big)\Big] = v(y) + v'(y)\varepsilon + \frac{1}{2}v''(y)\varepsilon^2 + o(\varepsilon^2)$$

$$= v(y + \varepsilon) + o(\varepsilon^2)$$

which follows from Lemma 4 in [14] and the identity (90) in Theorem 4.  □

Given Lemma 13, it is easy to complete the proof of Theorem 5. Indeed, using the notation of the lemma and defining the deterministic function

$$f(y) \triangleq \Big\|\eta(y) - \eta\Big(\frac{b}{2}\Big)\Big\|_{\mathbf{L}^2}, \qquad y < b,$$

we deduce from (107) that $f'(y) = 0$. As $f(\frac{b}{2}) = 0$, we have that $f(y) = 0$, $y < b$, and the result follows.  □

We recall that a financial model is called *complete* if it is arbitrage-free and any bounded, nonnegative contingent claim is replicable. It is well known that the model is complete if and only if the family of equivalent probability measures contains only one element $\mathbb{Q}$. Note that in this case, the density process of $\mathbb{Q}$ with respect to $\mathbb{P}$ is the greatest element of $\mathcal{Y}(1)$ in the sense that

$$(108) \qquad \frac{d\mathbb{Q}}{d\mathbb{P}} \geq Y_T \qquad \text{for any } Y \in \mathcal{Y}(1).$$

It is easy to show that if the financial model is complete, then the risk-tolerance wealth process is well defined for *any* utility function $U$ such that conditions (19) and (33) hold true. The next theorem provides a complete description of the general (incomplete) financial models sharing this property.

To state the result, we denote by $\succeq_2$ the *second-order stochastic dominance relation* between random variables. Recall that if $f$ and $g$ are *nonnegative* random variables on $(\Omega, \mathcal{F}, \mathbb{P})$, then $f \succeq_2 g$ if

$$\int_0^t \mathbb{P}(f \geq x)\, dx \geq \int_0^t \mathbb{P}(g \geq x)\, dx, \qquad t \geq 0.$$

It is well known that $f \succeq_2 g$ if and only if

$$\mathbb{E}[\phi(f)] \leq \mathbb{E}[\phi(g)]$$

for any function $\phi = \phi(x)$ on $[0, \infty)$ that is convex and decreasing and such that the expected values above are well defined.

THEOREM 6. *Assume* (11). *The following assertions are then equivalent:*

1. *The risk-tolerance wealth process* $R(x) = R(x|U)$ *is well defined for any initial wealth* $x > 0$ *and any utility function* $U$ *such that conditions* (19) *and* (33) *are satisfied.*

2. *The set* $\mathcal{Y}(1)$ *contains a greatest element* $\widehat{Y}$ *in the sense that*

$$(109) \qquad \widehat{Y}_T \succeq_2 Y_T \qquad \text{for any } Y \in \mathcal{Y}(1).$$

REMARK 3. If either of the equivalent assertions in Theorem 6 holds true, then the largest element of $\mathcal{Y}(1)$ is the density process of an equivalent *supermartingale* measure $\widehat{\mathbb{Q}}$, where the term "supermartingale measure" means that any $X \in \mathcal{X}(1)$ is a supermartingale under $\widehat{\mathbb{Q}}$. This follows from the definition of the set $\mathcal{Y}(1)$, condition (11) and the well-known fact that for nonnegative random variables $f$ and $g$, the relation $f \succeq_2 g$ implies that $\mathbb{E}[f] \geq \mathbb{E}[g]$. Note also that if the price process $S$ is sigma-bounded, then $\widehat{\mathbb{Q}}$ is, in fact, an equivalent martingale measure, that is, $\widehat{\mathbb{Q}} \in \mathcal{Q}$.

PROOF OF THEOREM 6. $2 \Longrightarrow 1$: Condition (109) implies that for any utility function $U$ satisfying the conditions of item 1, the dual value function is given by

$$v(y) = \mathbb{E}[V(y\widehat{Y}_T)],$$

where $V$ is the convex conjugate of $U$. The result now follows from Theorem 5.

$1 \Longrightarrow 2$: Denote by $\mathcal{V}$ the class of functions $V$ that are conjugate to the utility functions $U$ satisfying (19) and (33). We have that $V(y) \triangleq -\sqrt{y} \in \mathcal{V}$, as it is the convex conjugate of the function $U(x) \triangleq -1/4x$. Denote by $\widehat{Y}$ the element of $\mathcal{Y}(1)$ such that, for any $Y \in \mathcal{Y}(1)$, we have

$$(110) \qquad -\infty < -\mathbb{E}[\sqrt{\widehat{Y}_T}] \leq -\mathbb{E}[\sqrt{Y_T}].$$

To prove (109), it is sufficient to show that

$$(111) \qquad \mathbb{E}[\phi(\widehat{Y}_T)] \leq \mathbb{E}[\phi(Y_T)]$$

for any $Y \in \mathcal{Y}(1)$ and any bounded, decreasing, convex and two-times continuously differentiable function $\phi$ on $[0, \infty)$ such that $\phi(y) = 0$, $y \geq y_0$, for some $y_0 > 0$.

It is easy to see that if we select such a function $\phi$, then for any $\varepsilon > 0$,

$$\psi(y) \triangleq \phi(y) - \varepsilon\sqrt{y} \in \mathcal{V}.$$

From Theorem 5, we deduce the existence of $\widetilde{Y} \in \mathcal{Y}(1)$ such that for any $y > 0$ and $Y \in \mathcal{Y}(1)$,

$$\mathbb{E}[\psi(y\widetilde{Y}_T)] \leq \mathbb{E}[\psi(yY_T)].$$

Since $\phi$ is bounded, dividing this inequality by $\sqrt{y}$ and letting $y \to \infty$, we obtain

$$-\mathbb{E}[\sqrt{\widetilde{Y}_T}] \leq -\mathbb{E}[\sqrt{Y_T}], \qquad Y \in \mathcal{Y}(1)$$

and, hence, $\widetilde{Y} = \widehat{Y}$. Letting $\varepsilon \to 0$, we deduce (111) and thereby conclude the proof of the theorem. $\quad\square$

A natural question complementary to the one studied in Theorem 6 is to identify the class of utility functions $U$ such that the risk-tolerance wealth process is well defined for *any* finite financial model, that is, for any financial model where $|\Omega| < \infty$. In Theorem 7 below, we give this characterization for a class of utility functions more general than those studied thus far, namely for utility functions $U$ defined on $(a, \infty)$, where $-\infty \leq a < \infty$. As before, we assume that $U$ is strictly increasing, strictly convex, two-times differentiable on $(a, \infty)$ and satisfies the Inada conditions

$$(112) \qquad U'(a) = \lim_{x \to a} U'(x) = \infty, \qquad U'(\infty) = \lim_{x \to \infty} U'(x) = 0.$$

The next theorem says that the *only* utility functions allowing the existence of risk-tolerance wealth processes for an arbitrary financial model are power utilities, when $a > -\infty$, and exponential utilities, when $a = -\infty$.

THEOREM 7. *Let $-\infty \leq a < \infty$ and consider a utility function $U : (a, \infty) \to \mathbf{R}$ which is strictly increasing, strictly convex, two-times differentiable and which satisfies the Inada conditions* (112). *Then the following assertions are equivalent:*

1. *For any initial wealth $x > 0$ and any finite financial model that satisfies* (11), *the risk-tolerance wealth process $R(x)$ is well defined.*
2. *The utility function $U$ is such that*

   (a) *If $a = -\infty$,*

   $$U(x) = -c_1 e^{-\gamma x} + c_2, \qquad -\infty < x < \infty,$$

   *for some constants $c_1, \gamma > 0$ and $c_2 \in \mathbf{R}$,*

(b) *If $a > -\infty$, then either*

$$U(x) = c_1 \frac{(x-a)^p - 1}{p} + c_2, \qquad a < x < \infty,$$

*for some $c_1 > 0$, $p < 1$, $p \neq 0$ and $c_2 \in \mathbf{R}$, or*

$$U(x) = c_1 \ln(x-a) + c_2, \qquad a < x < \infty,$$

*for some $c_1 > 0$ and $c_2 \in \mathbf{R}$.*

PROOF.  Note first that both cases considered in item 2 have the same equivalent description in terms of the risk-tolerance coefficient of the utility function $U$:

$$(113) \qquad t(x) := -\frac{U'(x)}{U''(x)} = cx + d, \qquad x > a,$$

for some $c, d \in \mathbf{R}$.

$1 \Longrightarrow 2$: If the risk-tolerance coefficient $t(x)$ satisfies (113), it is obvious that the random variable $t(X_T(x))$ can be replicated and the replicating process is $cX(x) + d$.

$2 \Longrightarrow 1$: We prove this by contradiction. Assume that the risk-tolerance function $t(x)$ is not linear. Then there exist $0 < x_1 < x_2 < x_3 < \infty$ such that

$$(114) \qquad \frac{t(x_2) - t(x_1)}{x_2 - x_1} \neq \frac{t(x_3) - t(x_2)}{x_3 - x_2}.$$

Let us choose a one-period model with five possible outcomes, as follows:

$$S_0 = 1, \qquad S_1(\omega_1) = x_1, \qquad S_1(\omega_2) = x_2,$$

$$S_1(\omega_3) = x_3, \qquad S_1(\omega_4) = \tfrac{1}{2}, \qquad S_1(\omega_5) = 2.$$

We also choose the probabilities $p_1, p_2, p_3, p_4$ and $p_5$ of the outcomes such that

$$\mathbb{E}[U'(S_1)(S_1 - S_0)] = 0.$$

From the above equation, we conclude that the optimal strategy starting with $x = 1$ is to buy and hold one share of stock, that is, $X(1) = S$. As, by assumption, the random variable $t(S_1)$ can be replicated, there are $x$ and $\Delta$ in $\mathbf{R}$ such that

$$x + \Delta(S_1 - S_0) = t(S_1).$$

Listing the meaning of the above equation for the first three possible outcomes, we obtain

$$x + \Delta(x_1 - 1) = t(x_1),$$

$$x + \Delta(x_2 - 1) = t(x_2),$$

$$x + \Delta(x_3 - 1) = t(x_3),$$

which contradicts (114). The proof is thus complete.  □

**6. Qualitative properties of marginal utility-based prices.** As before, we fix an initial capital $x > 0$ and denote $y := u'(x)$. We assume that the conditions of Theorem 1 are satisfied. We start by computing $p'(x)$ and $D(x)$ under the condition that the risk-tolerance wealth process $R(x)$ is well defined.

By Theorem 4, $R(x)Y(y)$ is a martingale and, hence, we can define the probability measure $\widetilde{\mathbb{R}}(x)$ such that

$$(115) \qquad \frac{d\widetilde{\mathbb{R}}(x)}{d\mathbb{P}} = \frac{R_T(x)Y_T(y)}{R_0(x)y}.$$

Choosing $R(x)/R_0(x)$ as a numéraire, we denote by

$$(116) \qquad S^{R(x)} \triangleq \left( \frac{R_0(x)}{R(x)}, \frac{R_0(x)S}{R(x)} \right)$$

the price process of the traded securities and by

$$(117) \qquad \widetilde{g}_i(x) = \frac{f_i}{R_T(x)}R_0(x) = g_i(x)\frac{X_T(x)R_0(x)}{xR_T(x)}, \qquad 1 \le i \le m,$$

the payoffs of the contingent claims. Further, let

$$(118) \qquad \widetilde{\mathcal{M}}^2(x) \triangleq \left\{ M \in \mathbf{H}_0^2(\widetilde{\mathbb{R}}(x)) : M = \int H \, dS^{R(x)} \right\},$$

and $\widetilde{\mathcal{N}}^2(y)$ be the orthogonal complement to $\widetilde{\mathcal{M}}^2(x)$ in $\mathbf{H}_0^2(\widetilde{\mathbb{R}}(x))$. As we shall see in the proof of Theorem 8 below, the process

$$(119) \qquad \widetilde{P}_t(x) \triangleq \mathbb{E}_{\widetilde{\mathbb{R}}(x)}[\widetilde{g}(x)|\mathcal{F}_t], \qquad 0 \le t \le T,$$

is an $m$-dimensional square integrable martingale under $\widetilde{\mathbb{R}}(x)$. Noting that

$$(120) \qquad \widetilde{P}_0(x) = \mathbb{E}_{\widetilde{\mathbb{R}}(x)}[\widetilde{g}(x)] = \frac{1}{y}\mathbb{E}[Y_T(y)f] = p(x),$$

we denote by

$$(121) \qquad \widetilde{P}(x) = p(x) - \widetilde{M}(x) + \widetilde{N}(y)$$

the Kunita–Watanabe decomposition of $\widetilde{P}(x)$, where $\widetilde{M}^i(x)$ belongs to $\widetilde{\mathcal{M}}^2(x)$ and $\widetilde{N}^i(y)$ is an element of $\widetilde{\mathcal{N}}^2(y)$ [recall that $y = u'(x)$], $1 \le i \le m$.

THEOREM 8. *Let $x > 0$ and denote $y = u'(x)$. Assume that the conditions of Theorem 1 are satisfied and also that the risk-tolerance wealth process $R(x)$ exists. Then the marginal utility-based prices satisfy* (28) *and the corresponding sensitivity parameters $p'(x)$ and $D(x)$ equal*

$$(122) \qquad p'(x) = 0,$$

$$(123) \qquad D^{ij}(x) = \frac{u''(x)}{u'(x)}\mathbb{E}_{\widetilde{\mathbb{R}}(x)}[\widetilde{N}_T^i(y)\widetilde{N}_T^j(y)], \qquad 1 \le i, j \le m,$$

*where $\widetilde{N}(y)$ is defined by* (121). *The matrix $D(x)$ is symmetric negative semidefinite and if, in addition, Assumption* 1 *holds true, then it has full rank.*

PROOF.    The proof follows rather straightforwardly from Theorems 3 and 4. We start by showing that

$$\|\widetilde{g}(x)\| \in \mathbf{L}^2(\widetilde{\mathbb{R}}(x)) \tag{124}$$

and, hence, that $\widetilde{P}(x)$ is a square integrable martingale under $\widetilde{\mathbb{R}}(x)$.

By Assumption 2, there exist constants $c_1 > 0$ and $c_2 < \infty$ such that

$$c_1 R(x) \le X(x) \le c_2 R(x).$$

Assumption 3 and the above inequality imply that the process $S^{R(x)}$ defined in (116) is *sigma-bounded*. Furthermore, the sets $\mathcal{M}^2(x)$, $\widetilde{\mathcal{M}}^2(x)$ defined in (41) and (118) are related by

$$M \in \mathcal{M}^2(x) \quad \Longleftrightarrow \quad M \frac{X(x)R_0(x)}{R(x)x} \in \widetilde{\mathcal{M}}^2(x). \tag{125}$$

Assumption 4 now implies the existence of a constant $\widetilde{c} > 0$ and an $\widetilde{M} \in \widetilde{\mathcal{M}}^2(x)$ such that

$$\|\widetilde{g}(x)\| \le \widetilde{c} + \widetilde{M}_T$$

and (124) follows.

According to (82) and Theorem 4,

$$\widetilde{p}(x) = \mathbb{E}[Y'_T(y)f] = \frac{1}{y}\mathbb{E}[Y_T(y)f] = p(x).$$

The equality $p'(x) = 0$ now follows from (83). To prove the representation (123) of $D(x)$, we note that as

$$\frac{d\mathbb{R}(x)}{d\widetilde{\mathbb{R}}(x)} = \frac{X_T(x)R_0(x)}{R_T(x)x} = A(X_T(x))\frac{R_0(x)}{x},$$

the optimization problem (43) takes the form

$$a_i(x) \triangleq \inf_{M \in \mathcal{M}^2(x)} \mathbb{E}_{\mathbb{R}(x)}[A(X_T(x))(g_i(x) + M_T)^2]$$

$$= \frac{x}{R_0(x)} \inf_{\widetilde{M} \in \widetilde{\mathcal{M}}^2(x)} \mathbb{E}_{\widetilde{\mathbb{R}}(x)}[(\widetilde{g}_i(x) + \widetilde{M}_T)^2]$$

and, hence, its solution $M^i(x)$ is related to the process $\widetilde{M}^i(x)$ appearing in the Kunita–Watanabe decomposition (121) by

$$M^i(x) = \frac{R(x)x}{X(x)R_0(x)}\widetilde{M}^i(x).$$

Now, (123) follows from the formulas given in Theorems 1, 3 and 4 by direct computation. Finally, we clearly have that $D(x)$ is symmetric and negative semidefinite. It is also easy to see that under Assumption 1, this matrix has full rank. $\square$

As the assertions of Theorem 8 show, the existence of the risk-tolerance wealth process $R(x)$ provides affirmative answers to all of the questions of a qualitative nature formulated in Section 1. The following result shows that the reverse implication also holds true, that is, that if one of these questions has an affirmative answer for *any* family of contingent claims $f$, then $R(x)$ is well defined:

THEOREM 9. *Let $x > 0$ and denote $y = u'(x)$. Assume* (11) *and* (19) *and also that Assumptions 2 and 3 hold true. Then the following conditions are equivalent:*

1. *The risk-tolerance wealth process $R(x)$ exists.*
2. *For any family $f$ of contingent claims that satisfies Assumption 4, the sensitivity vector $p'(x)$ equals zero.*
3. *For any family $f$ of contingent claims that satisfies Assumption 4, the sensitivity matrix $D(x)$ is symmetric.*
4. *For any family $f$ of contingent claims that satisfies Assumption 4, the sensitivity matrix $D(x)$ is negative semidefinite, that is,* (5) *holds true.*
5. *For any family $f$ of contingent claims that satisfies Assumptions 1 and 4, the sensitivity matrix $D(x)$ has full rank.*

PROOF. $1 \Longrightarrow 2$, 3, 4 and 5: This follows directly from Theorem 8.

2 or $3 \Longrightarrow 1$: It is clear from (83) that condition 2 implies $\widetilde{p}(x) = p(x)$, so

$$\mathbb{E}[Y_T'(y)f] = \frac{1}{y}\mathbb{E}[Y_T(y)f] \tag{126}$$

for any one-dimensional contingent claim $f$ satisfying Assumption 4. Consider now a two-dimensional contingent claim $(f, \frac{X_T(x)}{x})$. In this case, $p_2(x) = \widetilde{p}_2(x) = 1$, so, using relation (84), we conclude that condition 3 also implies (126). Relation (126) suffices to conclude that

$$Y'(y) = \frac{Y(y)}{y}.$$

Since the processes $X'(x)$ and $Y'(y)$ are related by

$$U''(X_T(x))X_T'(x) = u''(x)Y_T'(y),$$

we can conclude that

$$X_T'(x) = -\frac{u''(x)}{u'(x)}\left(-\frac{U'(X_T(x))}{U''(X_T(x))}\right) > 0.$$

This implies that the risk-tolerance process $R(x)$ exists and

$$R(x) = -\frac{u'(x)}{u''(x)}X'(x).$$

$4 \Longrightarrow 2$: We use the notation $p'(x|f)$ and $D(x|f)$ to indicate the dependence of the sensitivity parameters on the particular family $f$ of contingent claims. From definition (27) of the sensitivity parameters, we deduce that for any vector $z \in \mathbf{R}^m$,

$$D\left(x\Big|f + \frac{X_T(x)}{x}z\right) = p'(x|f)z^T + D(x|f),$$

where $z^T$ is the transpose of $z$. Note that the random endowment $f + \frac{X_T(x)}{x}z$ satisfies Assumption 4. Since for any $z$, the matrix $D(x|f + \frac{X_T(x)}{x}z)$ is negative semidefinite, the last equality implies that $p'(x|f) = 0$.

$5 \Longrightarrow 1$: We choose the one-dimensional claim

$$(127) \qquad f = -\frac{U'(X_T(x))}{U''(X_T(x))} = -\frac{Y_T(y)}{U''(X_T(x))}.$$

We see that

$$x\frac{f}{X_T(x)} = x\frac{1}{A(X_T(x))},$$

so the solution of minimization problem (43) is given by $M(x) = 0$. According to (44) and Theorem 1, we obtain

$$(128) \qquad G_{11}(x) = \mathbb{E}[U''(X_T(x))(f)^2] = \mathbb{E}\left[\frac{(Y_T(y))^2}{U''(X_T(x))}\right].$$

We can also compute

$$(129) \qquad p(x) = \mathbb{E}\left[\frac{Y_T(y)}{y}f\right] = -\frac{1}{y}\mathbb{E}\left[\frac{(Y_T(y))^2}{U''(X_T(x))}\right]$$

and

$$\widetilde{p}(x) = \mathbb{E}\left[Y'_T(y)\frac{Y_T(y)}{-U''(X_T(x))}\right].$$

Since $U''(X_T(x))X'_T(x) = u''(x)Y'_T(y)$ and $X'(x)Y(y)$ is a uniformly integrable martingale, we obtain

$$(130) \qquad \widetilde{p}(x) = -\frac{y}{u''(x)}.$$

According to (84),

$$D(x) = \frac{1}{y}[G_{11}(x) - u''(x)p(x)\widetilde{p}(x)].$$

Using (128), (129) and (130), we obtain that $D(x) = 0$. This means that the contingent claim $f$ defined in (127) is replicable and, hence, the risk-tolerance wealth process exists. $\quad\square$

We conclude the section with an important remark:

Remark 4. The two key elements in the definition of marginal utility-based prices are the *financial model* (which consists of the price processes of traded securities defined on the filtered probability space) and the *utility function $U$* (which describes the preferences of the investor). By combining the results of Theorems 9, 6 and 7, we deduce that:

1. If the financial model is *fixed* and we want the qualitative properties (3)–(6) to be satisfied for *any* utility function, then the financial model should satisfy assertion 2 of Theorem 6, that is, the dual set $\mathcal{Y}(1)$ should have a greatest element from the point of view of second-order stochastic dominance.
2. If, on the other hand, the utility function $U$ is *fixed* and we want the qualitative properties (3)–(6) to be satisfied for *any* (finite) financial model, then the utility function should satisfy assertion 2 of Theorem 7, that is, $U$ should be a power utility (or an exponential utility if it is defined on the whole real line).

**7. An example.** In order to illustrate the general theory, we shall consider a specific example. We consider a European option written on a *nontraded, but observable* asset $Q$. We assume that the investor can trade an asset $S$ correlated to $Q$. A similar pricing problem was studied in [7] and later in [6] for power and exponential utilities.

Let $W = (W_t)_{0 \leq t \leq T}$ and $B = (B_t)_{0 \leq t \leq T}$ be two independent Brownian motions on a filtered probability space $(\Omega, \mathbb{P}, (\mathcal{F}_t)_{0 \leq t \leq T}, \mathcal{F})$, where the filtration is generated by $W$ and $B$. The evolution of the *nontraded* asset $Q$ is given by

$$(131) \qquad dQ_t = Q_t(\nu \, dt + \eta(\rho \, dW_t + \sqrt{1 - \rho^2} \, dB_t))$$

and the *traded* asset $S$ evolves according to

$$(132) \qquad dS_t = S_t(\mu \, dt + \sigma \, dW_t).$$

Here $\nu \in \mathbf{R}$, $\mu \in \mathbf{R}$, $\eta > 0$, $\sigma > 0$ and $0 < \rho < 1$ are constants. The money market pays zero interest rate.

Consider an economic agent starting with initial wealth $x > 0$ who can trade only in $S$. As before, we assume that the agent has a utility function $U$ satisfying Assumption 2. The agent is pricing a contingent claim with payoff $f = h(Q_T)$, where $h = h(x)$ is a bounded function. Of course, this covers the case of a European put written on $Q$.

We start with the verification of the conditions of Theorem 8. First, we show that the model satisfies the second condition of Theorem 6, where $\widehat{Y}$ is given by

$$d\widehat{Y}_t = -\frac{\mu}{\sigma}\widehat{Y}_t\,dW_t, \qquad \widehat{Y}_0 = 1.$$

Indeed, if $Y \in \mathcal{Y}(1)$, then using standard arguments based on the Doob–Meyer decomposition, the integral representation theorem and the computations behind Girsanov's theorem, we deduce that

$$Y = AN\widehat{Y},$$

where $A$ is a decreasing process, $A_0 = 1$ and $N$ satisfies

$$dN_t = N_t\alpha_t\,dB_t, \qquad N_0 = 1,$$

for some adapted process $\alpha$. Using the independence of $B$ and $W$, we deduce that $N$ is a nonnegative local martingale (hence, a supermartingale) under the enlarged filtration $(\mathcal{F}_t \vee \mathcal{F}_T^W)$. It follows that

$$\mathbb{E}[Y_T|\mathcal{F}_T^W] = \mathbb{E}[A_T N_T|\mathcal{F}_T^W]\widehat{Y}_T \leq \mathbb{E}[N_T|\mathcal{F}_T^W]\widehat{Y}_T \leq \widehat{Y}_T.$$

From Jensen's inequality, we deduce that for any bounded, decreasing and convex $\phi = \phi(x)$,

$$\mathbb{E}[\phi(Y_T)] \geq \mathbb{E}[\phi(\mathbb{E}[Y_T|\mathcal{F}_T^W])] \geq \mathbb{E}[\phi(\widehat{Y}_T)].$$

Hence, $\widehat{Y}_T \succeq_2 Y_T$, proving the second assertion of Theorem 6. This theorem then implies the existence of the risk-tolerance wealth process $R(x)$.

Assumption 3 in Theorem 8 holds true trivially, as $S$ is continuous. To verify the finiteness of $u = u(x)$ and Assumption 4, we recall (see Lemma 9 and the remark following it) that by Assumption 2, for any $0 < a < 1$, there are $0 < b_1 < b_2 < \infty$ such that

$$b_1 V''(y) \leq V''(ay) \leq b_2 V''(y),$$

where $V$ is the conjugate function to $U$. It follows that $V$, $V'$ and $V''$ are all bounded from below and above by power functions. Hence, for any $y > 0$,

$$(133) \qquad \mathbb{E}[V(y\widehat{Y}_T)] < \infty \quad \text{and} \quad \mathbb{E}\left[\frac{1}{V''(y\widehat{Y}_T)}\right] < \infty.$$

The first inequality in (133) readily implies the finiteness of the dual value function $v = v(y)$ and, hence, also the finiteness of $u = u(x)$. From the second inequality in (133) and Assumption 2 [see (34)], we deduce that

$$
\begin{aligned}
\mathbb{E}_{\mathbb{R}(x)}\left[\frac{x}{X_T^2(x)}\right] &= \mathbb{E}\left[\frac{\widehat{Y}_T}{X_T(x)}\right] = \mathbb{E}\left[\frac{\widehat{Y}_T}{-V'(y\widehat{Y}_T)}\right] \\
&\leq \frac{1}{yc_1}\mathbb{E}\left[\frac{1}{V''(y\widehat{Y}_T)}\right] < \infty,
\end{aligned}
\tag{134}
$$

where $X(x)$ is the optimal wealth process, $y = u'(x)$ and $\mathbb{R}(x)$ is the probability measure with density process $X(x)\widehat{Y}/x$. As $\widehat{Y}$ is a martingale under $\mathbb{P}$, it follows that $\frac{1}{X(x)}$ is a martingale under $\mathbb{R}(x)$. Then (134) implies that $\frac{1}{X(x)}$ is, in fact, a *square integrable martingale* under $\mathbb{R}(x)$ and Assumption 4 follows from the boundedness of the payoff $f = h(Q_T)$. This completes the verification of the conditions of Theorem 8.

To facilitate future computations, we introduce the process

$$
\widetilde{Q}_t = e^{-\kappa t}Q_t, \qquad 0 \leq t \leq T,
$$

where

$$
\kappa = \nu - \frac{\mu}{\sigma}\rho\eta.
$$

The benefit of this transformation comes from the fact that $\widetilde{Q}$ is a martingale (along with $S$) under the *minimal* martingale measure $\widehat{\mathbb{Q}}$ with density process $\widehat{Y}$. If $\widetilde{Q}$ were traded, then the contingent claim $f = h(e^{\kappa T}\widetilde{Q}_T)$ could be hedged by a bounded wealth process $P$:

$$
P_t = \mathbb{E}_{\widehat{\mathbb{Q}}}[h(e^{\kappa T}\widetilde{Q}_T)|\mathcal{F}_t] = p(x) + \int_0^t \Delta(\widetilde{Q}_s, s)\frac{d\widetilde{Q}_s}{\widetilde{Q}_s},
$$

where $\Delta = \Delta(x, t)$ is the amount of money invested in $\widetilde{Q}$ at $t$ if $\widetilde{Q}_t = x$.

Following the path outlined in Theorem 8, we denote by $\widetilde{P}(x) = PR_0(x)/R(x)$ the price process of the option discounted by the normalized risk-tolerance process and by $\widetilde{\mathbb{R}}(x)$ the probability measure given by

$$
\frac{d\widetilde{\mathbb{R}}(x)}{d\mathbb{P}} = \frac{R_T(x)\widehat{Y}_T}{R_0(x)}.
$$

Since $\frac{d\widetilde{\mathbb{R}}(x)}{d\mathbb{P}}$ is $\mathcal{F}_T^W$-measurable, we can decompose the process $W$ into $\widetilde{W}$ plus a "drift" such that $\widetilde{W}$ and $B$ are independent Brownian motions under the measure $\widetilde{\mathbb{R}}(x)$. The two-dimensional stock process $S^{R(x)}$ defined in (116) becomes a stochastic integral with respect to $\widetilde{W}$ only. Straightforward computations (identifying the "$d\widetilde{W}_t$" and the "$dB_t$"-parts) show that

$$
d\widetilde{P}_t(x) = H_t\, dS_t^{R(x)} + \eta\sqrt{1-\rho^2}\Delta(\widetilde{Q}_t, t)\frac{R_0(x)}{R_t(x)}\, dB_t
$$

for some (two-dimensional) integrand $H$. This is the desired Kunita–Watanabe decomposition (121). Using the fact that $R_0(x) = -u'(x)/u''(x)$, we obtain from Theorem 8 that $p'(x) = 0$ and

$$(135) \qquad D(x) = \eta^2(1-\rho^2)\frac{u'(x)}{u''(x)}\mathbb{E}_{\widetilde{\mathbb{R}}(x)}\left[\int_0^T \left(\frac{\Delta(\widetilde{Q}_t,t)}{R_t}\right)^2 dt\right].$$

For the case of power utilities, this formula has been obtained in [6], (4.19).

## APPENDIX A: SECOND-ORDER EXPANSIONS FOR RESERVATION PRICES AND CERTAINTY EQUIVALENCE VALUES

For the convenience of future references, we present here the asymptotic analysis for the quantities $b(x,q)$ and $c(x,q)$ defined as the solutions to the following equations:

$$(136) \qquad\qquad u(x) = u(x - b(x,q), q),$$

$$(137) \qquad\qquad u(x,q) = u(x + c(x,q)).$$

Here $(x,q) \in \mathcal{K}$, $x > 0$ and the value functions $u(x)$ and $u(x,q)$ are defined in (18) and (17). In the literature, $b(x,q)$ is often referred to as the *reservation price* of the contingent claim $\langle q, f \rangle$ and the quantity $x + c(x,q)$ is called the *certainty equivalence value* of the portfolio $(x,q)$.

THEOREM A.1. *Fix $x > 0$ and assume the conditions and the notation of Theorem 3. Then the quantities $b(x,q)$ and $c(x,q)$ defined in (136) and (137) have the following second-order expansions around $(x,0)$:*

$$(138) \begin{aligned} b(x + \Delta x, q) &= \langle p(x), q \rangle + \tfrac{1}{2}\begin{pmatrix} \Delta x & q \end{pmatrix} B(x) \begin{pmatrix} \Delta x \\ q \end{pmatrix} \\ &\quad + o(|\Delta x|^2 + \|q\|^2), \end{aligned}$$

$$(139) \begin{aligned} c(x + \Delta x, q) &= \langle p(x), q \rangle + \tfrac{1}{2}\begin{pmatrix} \Delta x & q \end{pmatrix} C(x) \begin{pmatrix} \Delta x \\ q \end{pmatrix} \\ &\quad + o(|\Delta x|^2 + \|q\|^2), \end{aligned}$$

*where $B(x)$ and $C(x)$ are $(m+1) \times (m+1)$ symmetric matrices given by the formulas (for $1 \leq i, j \leq m$):*

$$(140) \begin{aligned} B_{00}(x) &= C_{00}(x) = 0, \\ B_{0i}(x) &= C_{0i}(x) = p'_i(x), \\ B_{ij}(x) &= D_{ij}(x) - p'_i(x)p_j(x), \\ C_{ij}(x) &= D_{ij}(x) + p_i(x)p'_j(x). \end{aligned}$$

*In addition, if the risk-tolerance wealth process $R(x)$ is well defined, then for $1 \le i, j \le m$,*

$$
\begin{aligned}
B_{0i}(x) &= C_{0i}(x) = 0, \\
B_{ij}(x) &= C_{ij}(x) = D_{ij}(x).
\end{aligned}
\tag{141}
$$

PROOF. Relations (138) and (139) are rather straightforward corollaries of the second-order approximation (29) for $u(x, q)$ if we take into account the notation used in Theorem 3. Indeed, considering, for example, the case of the "liquidation" price $c(x, q)$ and denoting

$$
\widetilde{c}(x + \Delta x, q) \triangleq \langle p(x), q \rangle + \tfrac{1}{2} \begin{pmatrix} \Delta x & q \end{pmatrix} C(x) \begin{pmatrix} \Delta x \\ q \end{pmatrix},
$$

where the matrix $C(x)$ is given by (140), we deduce from (29) by direct computation that

$$
u(x + \Delta x + \widetilde{c}(x + \Delta x, q)) - u(x + \Delta x, q) = o(|\Delta x|^2 + \|q\|^2).
$$

As $u'(x) > 0$, there is a constant $K > 0$ such that for sufficiently small $\Delta x$ and $q$,

$$
\begin{aligned}
K|\widetilde{c}&(x + \Delta x, q) - c(x + \Delta x, q)| \\
&\le |u(x + \Delta x + \widetilde{c}(x + \Delta x, q)) - u(x + \Delta x + c(x + \Delta x, q))| \\
&= |u(x + \Delta x + \widetilde{c}(x + \Delta x, q)) - u(x + \Delta x, q)| = o(|\Delta x|^2 + \|q\|^2),
\end{aligned}
$$

which implies (139). Finally, relations (141) follow from Theorem 8. □

## APPENDIX B: ON THE ATTAINABILITY OF BOUNDED CONTINGENT CLAIMS

We consider the setting of Section 2 and denote by $\mathcal{R} = \mathcal{R}(S)$ the family of equivalent probability measures $\mathbb{R}$ such that any $X \in \mathcal{X}(1)$ is a *supermartingale* under $\mathbb{R}$. Note that if $S$ is locally bounded, then $\mathcal{R}$ coincides with the family $\mathcal{Q} = \mathcal{Q}(S)$ of equivalent local martingale measures for $S$ and that, in general, $\mathcal{Q} \subset \mathcal{R}$. To simplify notation we assume that

$$
\mathbb{P} \in \mathcal{R}
\tag{142}
$$

and denote by $\mathcal{R}^\infty$ the elements of $\mathcal{R}$ that have bounded densities with respect to $\mathbb{P}$, that is,

$$
\mathcal{R}^\infty = \left\{ \mathbb{R} \in \mathcal{R} : \frac{d\mathbb{R}}{d\mathbb{P}} \in \mathbf{L}^\infty \right\}.
$$

We also denote by $\mathcal{R}^0$ the set

$$
\mathcal{R}^0 = \{ \mathbb{R} \in \mathcal{R} : \mathbb{P} \pm \varepsilon(\mathbb{R} - \mathbb{P}) \in \mathcal{R} \text{ for some } \varepsilon > 0 \}.
$$

It is obvious that

$$\mathcal{R}^0 \subset \mathcal{R}^\infty \subset \mathcal{R}.$$

The following result is a generalization of the well-known criteria of attainability of bounded contingent claims in incomplete markets established for the case when $\mathbb{P} \in \mathcal{Q}$; see Theorem 3.2 in [11].

**THEOREM B.1.** *Assume* (142) *and let* $\zeta \in \mathbf{L}^\infty$. *Then the contingent claim* $\zeta$ *is replicable if and only if*

$$(143) \qquad\qquad \mathbb{E}_\mathbb{R}[\zeta] = \mathbb{E}[\zeta] \qquad \forall \mathbb{R} \in \mathcal{R}^0.$$

The proof of Theorem B.1 relies on some auxiliary results. For a $d$-dimensional process $X$, we denote by $X^*$ the corresponding maximal process defined by

$$X_t^* \triangleq \sup_{0 \le s \le t} \|X_s\|, \qquad 0 \le t \le T.$$

**LEMMA B.1.** *Let* $X$ *be a* $d$-*dimensional special semimartingale. Then there exists a predictable, strictly positive, one-dimensional process* $\phi$ *such that*

$$\mathbb{E}\left[\left(\int \phi\, dX\right)_T^*\right] \le 1.$$

PROOF. The proof follows easily from the well-known fact that $X$ is a special semimartingale if and only if the maximal process $X^*$ is locally integrable. □

Proposition B.1 below is closely related to [14], Proposition 1. For a $d$-dimensional semimartingale $R$, denote

$$\mathcal{I}(R) = \left\{ X : X = \int H\, dR \text{ for some predictable } H \right\},$$

$$\mathcal{I}^{\mathrm{sp}}(R) = \{ X \in \mathcal{I}(R) : X \text{ is a special semimartingale} \}.$$

**PROPOSITION B.1.** *Let* $S$ *be a* $d$-*dimensional semimartingale. There exists a* $d$-*dimensional special semimartingale* $\widetilde{S} = (\widetilde{S}^i)_{1 \le i \le d}$ *such that* $\widetilde{S}^i \in \mathcal{I}^{\mathrm{sp}}(S)$ *for each* $1 \le i \le d$ *and*

$$\mathcal{I}^{\mathrm{sp}}(S) = \mathcal{I}^{\mathrm{sp}}(\widetilde{S}).$$

PROOF. We closely follow the proof of [14], Proposition 1. Let $\mathbb{Q}$ be an equivalent probability measure such that $S$ is a special semimartingale under $\mathbb{Q}$. We denote by

$$S = S_0 + M + A$$

the canonical decomposition of $S$ under $\mathbb{Q}$, where $M$ is a local martingale under $\mathbb{Q}$ and $A$ is a predictable process of bounded variation. Without any loss of generality (integrating, if necessary, by a one-dimensional strictly positive process), we can assume that

$$(144) \qquad M \in \mathbf{H}^1(\mathbb{Q}) \quad \text{and} \quad \mathbb{E}_{\mathbb{Q}}\left[\int_0^T \|dA_t\|\right] < \infty.$$

We denote by $C_t$ the one-dimensional compensator of $\|[M, M]\|^{1/2}$ under $\mathbb{Q}$ and define the measure $\mu$ on the predictable $\sigma$-field of $[0, T] \times \Omega$ by

$$d\mu(t, \omega) = (dC_t(\omega) + \|dA_t(\omega)\|) \, d\mathbb{Q}(\omega).$$

According to (144), the measure $\mu$ is finite, that is,

$$\mu([0, T] \times \Omega) < \infty.$$

Denote by $\mathcal{A}^{\mathrm{sp}}$ the set of predictable processes with values in the set of $d \times d$ symmetric and positive semidefinite matrices such that $X := \int A \, dS$ is well defined and $X$ is a special semimartingale under $\mathbb{P}$. We claim that there exists $\widetilde{A} \in \mathcal{A}^{\mathrm{sp}}$ such that

$$(145) \qquad \operatorname{rank}(\widetilde{A}) \geq \operatorname{rank}(A) \qquad \text{for each } A \in \mathcal{A}^{\mathrm{sp}},$$

where $\operatorname{rank}(A)$ denotes the rank of the matrix $A$. In order to prove this claim, we choose a sequence $(A^n)_{n \geq 1}$ in $\mathcal{A}^{\mathrm{sp}}$ such that

$$(146) \qquad \lim_{n \to \infty} \int_{[0,T] \times \Omega} \operatorname{rank}(A^n) \, d\mu = \sup_{A \in \mathcal{A}^{\mathrm{sp}}} \int_{[0,T] \times \Omega} \operatorname{rank}(A) \, d\mu.$$

According to Lemma B.1, $(A^n)_{n \geq 1}$ can be chosen such that

$$(147) \qquad \mathbb{E}\left[\left(\int A^n \, dS\right)_T^*\right] \leq 1, \qquad n \geq 1.$$

Of course, we can select $A^n$ so that, in addition, $\|A^n\| \leq 1$. In this case, we can define

$$\widetilde{A} \triangleq \sum_{n=1}^{\infty} \frac{1}{2^n} A^n.$$

Clearly,

$$\operatorname{rank}(\widetilde{A}) \geq \sup_{n \geq 1} \operatorname{rank}(A^n).$$

Moreover, by (147),

$$\mathbb{E}\left[\left(\int \widetilde{A} \, dS\right)_T^*\right] \leq \sum_{n=1}^{\infty} \frac{1}{2^n} \mathbb{E}\left[\left(\int A^n \, dS\right)_T^*\right] \leq 1.$$

Hence, $\widetilde{A} \in \mathcal{A}^{\mathrm{sp}}$ and satisfies (145).

Denote $\widetilde{S} \triangleq \int \widetilde{A} \, dS$. Let $H$ be a $d$-dimensional predictable process such that $\int H \, dS$ is a special semimartingale. To complete the proof, it remains to be shown that $H$ admits the representation

$$(148) \qquad\qquad H = \widetilde{A} G$$

for some predictable $d$-dimensional process $G$. Indeed, in this case,

$$\int H \, dS = \int \langle H, dS \rangle = \int \langle \widetilde{A} G, dS \rangle = \langle G, \widetilde{A} \, dS \rangle = \int G \, d\widetilde{S}.$$

To prove (148), we use the fact that $H$ can always be decomposed as

$$(149) \qquad\qquad H = \widetilde{A} G + F,$$

where $G$ and $F$ are predictable $d$-dimensional processes such that $F \in \ker(\widetilde{A})$ [$\ker(A)$ denotes the *kernel* of the matrix $A$]. Multiplying, if necessary, both sides of (149) by a strictly positive predictable process, we can assume that $\int F \, dS$ is a special semimartingale and $\|F\| \leq 1$. In this case, the matrix $B$ defined by

$$B^{ij} = F^i F^j, \qquad 1 \leq i, j \leq d,$$

belongs to $\mathcal{A}^{\mathrm{sp}}$. Hence, $\widetilde{A} + B \in \mathcal{A}^{\mathrm{sp}}$. However, as $F \in \ker(\widetilde{A})$,

$$\operatorname{rank}(\widetilde{A} + B) = \operatorname{rank}(\widetilde{A}) + I_{\{F \neq 0\}}$$

and the fact that $F$ equals zero follows from the maximality property (145) for $\widetilde{A}$. $\quad\square$

PROOF OF THEOREM B.1. If $\zeta$ is replicable, then (143) follows easily since any bounded $X \in \mathcal{X}(1)$ is a martingale under any $\mathbb{R} \in \mathcal{R}$. Below, we shall prove the reverse implication.

Since $\mathbb{P} \in \mathcal{R}$, it follows that each $X \in \mathcal{X}(1)$ is a positive supermartingale under $\mathbb{P}$, hence a special semimartingale. According to Proposition B.1, we can assume that $S$ is a special semimartingale under $\mathbb{P}$, having canonical decomposition

$$S = M + A,$$

where $M$ is a local martingale and $A$ is a predictable process of bounded variation, $A_0 = 0$. Denote by $\mathcal{Q}(M)$ the family of equivalent local martingale measures for $M$. We claim that

$$(150) \qquad\qquad \mathcal{Q}(M) \subset \mathcal{R}.$$

Indeed, let $\int H\,dS$ be bounded from below. Since $\mathbb{P} \in \mathcal{R}$, this stochastic integral is a supermartingale and, as is well known, its Doob–Meyer decomposition is given by

$$\int H\,dS = \int H\,dM + \int H\,dA.$$

The stochastic integral $\int H\,dA$ being locally bounded, $\int H\,dM$ is locally bounded from below. It follows that $\int H\,dM$ is a local martingale under any $\mathbb{Q} \in \mathcal{Q}(M)$. As $\int H\,dA$ is decreasing, we deduce that $\int H\,dS$ is a supermartingale under any $\mathbb{Q} \in \mathcal{Q}(M)$ and (150) follows.

Denote now

$$\mathcal{Q}^0(M) = \left\{ \mathbb{Q} \in \mathcal{Q}(M) : \frac{d\mathbb{Q}}{d\mathbb{P}} \in \mathbf{L}^\infty \right\}.$$

Since $\mathbb{P} \pm \varepsilon(\mathbb{R} - \mathbb{P}) \in \mathcal{Q}(M)$ for any $\mathbb{R} \in \mathcal{Q}^0(M)$ and $\varepsilon > 0$ sufficiently small, using (150) can we deduce that $\mathcal{Q}^0(M) \subset \mathcal{R}^0$. Hence, by (143),

$$(151) \qquad \mathbb{E}_{\mathbb{Q}}[\zeta] = \mathbb{E}[\zeta] \qquad \forall \mathbb{Q} \in \mathcal{Q}^0(M).$$

By Theorem 1.1 in [11], condition (151) implies the existence of a constant $x$ and a bounded stochastic integral $\int H\,dM$ such that

$$\zeta = x + \int_0^T H_t\,dM_t.$$

Hence, the proof will be complete if we can show that $H$ is $S$-integrable and that

$$(152) \qquad \int H\,dM = \int H\,dS.$$

We claim that, in fact, (152) is valid for any predictable process $H$ such that $\int H\,dM$ is well defined and locally bounded. Indeed, if, in addition, $H$ is $A$-integrable, then the stochastic integral $\int H\,dS = \int H\,dM + \int H\,dA$ is well defined and locally bounded. It follows that $\int H\,dS$ is a local martingale and, hence, $\int H\,dA = 0$. In the general case, we fix $n \geq 0$ and denote $H^n = HI(|H| \leq n)$. Clearly, $\int H^n\,dM$ is locally bounded and $\int H^n\,dA$ is well defined. Hence, by what we have already proven,

$$\int H^n\,dA = \int HI(|H| \leq n)\,dA = 0$$

and the result follows because $n$ is an arbitrary positive number. $\square$

**Acknowledgments.** Part of the work was done while Mihai Sîrbu was a graduate student at Carnegie Mellon University.

DEPARTMENT OF MATHEMATICAL SCIENCES
CARNEGIE MELLON UNIVERSITY
5000 FORBES AVENUE
PITTSBURGH, PENNSYLVANIA 15213
USA
E-MAIL: kramkov@andrew.cmu.com

DEPARTMENT OF MATHEMATICS
COLUMBIA UNIVERSITY
2990 BROADWAY
NEW YORK, NEW YORK 10027
USA
E-MAIL: sirbu@math.columbia.edu